\documentclass[final]{siamart1116}

\usepackage{amsfonts,amssymb,amsmath,amsopn,float}
\usepackage{xifthen}
\usepackage{graphicx}
\usepackage{tikz}
\usepackage{pgfplots}
\usepgfplotslibrary{units}
\usepackage{epstopdf}
\usepackage{cite}
\usepackage{url}

\newcommand{\bbR}{\mathbb{R}}

\newcommand{\bbZ}{\mathbb{Z}}

\newcommand{\bzero}{\mathbf{0}}
\newcommand{\abs}[1]{\left\lvert #1 \right\rvert}
\newcommand{\norm}[1]{\left\lVert #1 \right\rVert}
\newcommand{\normX}[2]{{\left\lVert #1 \right\rVert}_{#2}}
\newcommand{\del}{\boldsymbol{\nabla}}

\newcommand{\leftcr}{\left\{}
\newcommand{\rightcr}{\right\}}
\newcommand{\bxi}{\boldsymbol{\xi}}
\newcommand{\betanew}{\boldsymbol{\eta}}

\newcommand{\dparder}[2]{\dfrac{\partial #1}{\partial #2}}

\newcommand{\dsecder}[2]{\dfrac{\partial^2 #1}{\partial #2^2}}

\newcommand{\Ltwo}[1]{%
\ifthenelse{\equal{#1}{}}{L^2}{L^2(#1)}%
}

\newcommand{\Ltwoz}[1]{%
\ifthenelse{\equal{#1}{}}{L^2_0}{L^2_0(#1)}%
}

\newcommand{\Ltwonorm}[2]{\left\lVert #1 \right\rVert_{\Ltwo{#2}}}

\newcommand{\Cone}[1]{%
\ifthenelse{\equal{#1}{}}{C^{1}}{C^{1}(#1)}%
}

\newcommand{\Conez}[1]{%
\ifthenelse{\equal{#1}{}}{C^{1}_{0}}{C^{1}_{0}(#1)}%
}

\newcommand{\Ctwo}[1]{%
\ifthenelse{\equal{#1}{}}{C^{2}}{C^2(#1)}%
}

\newcommand{\Ctwoz}[1]{%
\ifthenelse{\equal{#1}{}}{C^{2}_{0}}{C^{2}_{0}(#1)}%
}

\newcommand{\Cholder}[1]{%
\ifthenelse{\equal{#1}{}}{C^{0,\gamma}}{C^{0,\gamma}(#1)}%
}

\newcommand{\Cholderz}[1]{%
\ifthenelse{\equal{#1}{}}{C^{0,\gamma}_{0}}{C^{0,\gamma}_{0}(#1)}%
}

\newcommand{\Choldernorm}[2]{%
{\left\lVert #1 \right\rVert}_{\Cholder{#2}}%
}

\newcommand{\Ctwointime}[2]{C^{2}(#1;#2)}
\newcommand{\Cthreeintime}[2]{C^{3}(#1;#2)}

\newcommand{\perdrd}{D;\bbR^{d}}

\newcommand{\bolds}[1]{\boldsymbol{#1}}
\newcommand{\ba}{\bolds{a}}
\newcommand{\bb}{\bolds{b}}

\newcommand{\be}{\bolds{e}}

\newcommand{\br}{\bolds{r}}
\newcommand{\bs}{\bolds{s}}

\newcommand{\bu}{\bolds{u}}
\newcommand{\bv}{\bolds{v}}

\newcommand{\bx}{\bolds{x}}
\newcommand{\by}{\bolds{y}}
\newcommand{\bz}{\bolds{z}}

\newcommand{\bW}{\bolds{W}}

\newcommand{\buhat}{\hat{\bu}}
\newcommand{\bubar}{\bar{\bu}}
\newcommand{\butilde}{\tilde{\bu}}
\newcommand{\bvhat}{\hat{\bv}}
\newcommand{\bvbar}{\bar{\bv}}
\newcommand{\bvtilde}{\tilde{\bv}}

\def\perpotfx(#1,#2,#3){#2*#1/(1+#3*#1)}
\def\derperpotfx(#1,#2,#3){#2/(1+#3*#1) - #2*#3*#1/((1+#3*#1)^2)}
\def\dderperpotfx(#1,#2,#3){-2*#2*#3/((1+#3*#1)^3)}
\def\ddderperpotfx(#1,#2,#3){6*#2*#3*#3/((1+#3*#1)^4)}
\def\perpotFx(#1){\perpotfx(#1*#1,0.8,0.1)}
\def\derperpotFx(#1){2*#1*\derperpotfx(#1*#1,0.8,0.1)}
\def\dderperpotFx(#1){2(\derperpotfx(#1*#1,4.0,0.5) + 2*#1*#1*\dderperpotfx(#1*#1,4.0,0.5))}
\def\ddderperpotFx(#1){4*#1*(2*\dderperpotfx(#1*#1,0.8,0.1) + #1*#1*\ddderperpotfx(#1*#1,2.5,0.1))}

\newcommand{\TheTitle}{Numerical analysis of nonlocal fracture models in H\"{o}lder space} 
\newcommand{\TheAuthors}{Prashant K. Jha and Robert Lipton}

\title{{\TheTitle}\thanks{Published on SIAM Journal on Numerical Analysis on April 10, 2018.
\funding{This material is based upon work supported by the U. S. Army Research Laboratory and the U. S. Army Research Office under contract/grant number W911NF1610456.}}}

\author{
  Prashant K. Jha
  \thanks{Department of Mathematics, Louisiana State University, Baton Rouge, LA 
  (\email{prashant.j16o@gmail.com}).}
  \and
  Robert Lipton
  \thanks{Department of Mathematics, Louisiana State University, Baton Rouge, LA 
  (\email{lipton@math.lsu.edu}).}
}

\headers{\TheTitle}{\TheAuthors}

\ifpdf
  \DeclareGraphicsExtensions{.eps,.pdf,.png,.jpg}
\else
  \DeclareGraphicsExtensions{.eps}
\fi

\ifpdf
\hypersetup{
  pdftitle={\TheTitle},
  pdfauthor={\TheAuthors}
}
\fi

\begin{document}

\maketitle

\begin{abstract}
In this work, we calculate the convergence rate of the finite difference approximation for a class of nonlocal fracture models. We consider two point force interactions characterized by a double well potential. We show the existence of a evolving displacement field  in H\"{o}lder space with H\"{o}lder exponent $\gamma \in (0,1]$. 
The rate of convergence of the finite difference approximation depends on the factor $C_s h^\gamma/\epsilon^2$ where $\epsilon$ gives the length scale of nonlocal interaction, $h$ is the discretization length and $C_s$ is the maximum of H\"older norm of the solution and its second derivatives during the evolution. It is shown that the rate of convergence holds for both the  forward Euler scheme as well as general single step implicit schemes. A stability result is established for the semi-discrete approximation. The H\"older continuous evolutions are seen to converge to a brittle fracture evolution in the limit of vanishing nonlocality.
\end{abstract}

\begin{keywords}
  Nonlocal fracture models, peridynamics, cohesive dynamics, numerical analysis, finite difference approximation
\end{keywords}

\begin{AMS}
      34A34, 34B10, 74H55, 74S20
\end{AMS}

\section{Introduction}
Nonlocal formulations have been proposed to describe the evolution of deformations which exhibit loss of differentiability and continuity, see \cite{CMPer-Silling} and \cite{States}. These models are commonly referred to as peridynamic models. The main idea is to define the strain in terms of displacement differences and allow nonlocal interactions between material points. This generalization of strain allows for the participation of a larger class of deformations in the dynamics.  Numerical simulations based on peridynamic modeling exhibit formation and evolution of sharp interfaces associated with phase transformation and fracture  \cite{CMPer-Dayal}, \cite{CMPer-Silling4}, \cite{CMPer-Silling5}, \cite{CMPer-Silling7}, \cite{CMPer-Agwai}, \cite{CMPer-Du}, \cite{CMPer-Lipton2}, \cite{BobaruHu}, \cite{HaBobaru}, \cite{SillBob}, \cite{WeckAbe}, \cite{GerstleSauSilling}. A recent summary of the state of the art can be found in \cite{Handbook}.

In this work, we provide a numerical analysis for the class of nonlocal models introduced in \cite{CMPer-Lipton3}  and \cite{CMPer-Lipton}. These models are  defined by a double well two point potential. Here one potential well is centered at zero and associated with elastic response while the other well is at infinity and associated with surface energy. The rational for studying these models is that they are shown to be well posed over the class of square integrable non-smooth displacements and, in the limit of vanishing non-locality,  their dynamics recover features associated with sharp fracture propagation see, \cite{CMPer-Lipton3}  and \cite{CMPer-Lipton}.  The numerical simulation of prototypical fracture problems using this model is carried out in \cite{CMPer-Lipton2}.  In order to develop an $L^2$ approximation theory, we show the nonlocal evolution is well posed over a more regular space of functions. To include displacement fields which have no well-defined derivatives, we consider displacement fields in the H\"{o}lder space $\Cholder{}$ with H\"{o}lder exponent $\gamma$ taking any value in $(0,1]$. We show that a unique evolution exists in $\Cholder{}$ for  $\Cholder{}$  initial data and body force.  
The semi-discrete approximation to the H\"{o}lder continuous evolution is considered and it is shown that at any time its energy is bounded by the initial energy and the work done by the body force.
We develop an approximation theory for the forward Euler scheme and  show that these ideas can be easily extended to the backward Euler scheme as well other implicit one step time discretization schemes. 
It is  found that the discrete approximation converges to the exact solution in the $L^2$ norm uniformly over finite time intervals with the rate of convergence proportional to  $(C_t\Delta t + C_s h^\gamma/\epsilon^2)$, where $\Delta t$ is the size of time step, $h$ is the size of spatial mesh discretization, and $\epsilon$ is the length scale of nonlocal interaction relative to the size of the domain. The constant $C_t$ depends on the $L^2$ norm of the time derivatives of the solution,  $C_s$ depends on the H\"older norm of the solution and the Lipschitz constant of peridynamic force. 
We point out that the constants appearing in the convergence estimates with respect to $h$ can be dependent on the horizon and be large when $\epsilon$ is small.  This is discussed in  \autoref{s:finite difference} and an example is provided in \autoref{s:conclusions}. These results show that while errors can grow with each time step they can be controlled over finite times $t$ by suitable spatial temporal mesh refinement. We then apply the methods developed in \cite{CMPer-Lipton3}  and \cite{CMPer-Lipton}, to show that in the limit $\epsilon\rightarrow 0$, the H\"older continuous evolutions converge to a limiting sharp fracture evolution with bounded Griffiths fracture energy. Here the limit evolution is differentiable off the crack set and satisfies the linear elastic wave equation. 

In the language of nonlocal operators, the integral kernel associated with the nonlocal model studied here is Lipschitz continuous guaranteeing global stability of the finite difference approximation. This is in contrast to PDE based evolutions where stability can be conditional. 
In addition we examine local stability. Unfortunately the problem is nonlinear so we don't establish  CFL conditions but instead identify a mode of dynamic instability that can arise during the evolution. This type of instability is due to a radial perturbation of the solution and causes error to grow with each time step for the Euler scheme. For implicit schemes this perturbation can become unstable in parts of the computational domain where there is material softening, see \autoref{ss:proof localstab}. Of course stability conditions like the CFL conditions for linear nonlocal equations are of importance for  guidance in implementations. In the case of $d=1$, a CFL type condition is obtained for the finite difference and finite element approximation of the linear peridynamic equation, see \cite{CMPer-Guan}. Recent work develops a new simple CFL condition for one dimensional linearized   peridynamics in the absence of body forces \cite{CMPer-JhaLipton}. Related analysis for the linear peridynamic equation in one dimension is taken up in \cite{CMPer-Weckner} and \cite{CMPer-Bobaru}. The recent and related work  \cite{CMPer-Du1} and \cite{CMPer-Guan2} addresses numerical approximation for problems of nonlocal diffusion. 

There is now a large body of contemporary work addressing the numerical approximation of singular kernels with application to nonlocal diffusion, advection, and mechanics. Numerical formulations and convergence theory for nonlocal $p$-Laplacian formulations are developed in \cite{DeEllaGunzberger}, \cite{Nochetto1}. Numerical analysis of nonlocal steady state diffusion is presented in \cite{CMPer-Du2} and \cite{CMPer-Du3}, and \cite{CMPer-Chen}.  The use of fractional Sobolev spaces for nonlocal problems is investigated and developed in \cite{CMPer-Du1}. Quadrature approximations and stability conditions for linear peridynamics are analyzed in \cite{CMPer-Weckner} and  \cite{CMPer-Silling8}.  The interplay between nonlocal interaction length and grid refinement for linear peridynamic models is presented in \cite{CMPer-Bobaru}. Analysis of adaptive refinement and domain decomposition for linearized peridynamics are provided in \cite{AksoyluParks}, \cite{LindParks}, and \cite{AksMen}. This list is by no means complete and the literature on numerical methods and analysis continues to grow.

The paper is organized as follows. In \autoref{s:peridynamic model}, we describe the nonlocal model. In \autoref{ss:existence holder}, we state theorems which show Lipschitz continuity of the nonlocal force (\autoref{prop:lipschitz}) and the existence and uniqueness of an evolution over any finite time interval (\autoref{thm:existence over finite time domain}). In \autoref{s:finite difference}, we compute the convergence rate of the forward Euler scheme as well as implicit one step methods. We identify stability of the semi-discrete approximation with respect to the energy in \autoref{semidiscrete}. In \autoref{ss:proof localstab}, we identify local instabilities in the fully discrete evolution caused by suitable radial perturbations of the solution.  In \autoref{s:proof existence}, we give the proof of \autoref{prop:lipschitz}, \autoref{thm:local existence}, and \autoref{thm:existence over finite time domain}.  The convergence of H\"older continuous evolutions to sharp fracture evolutions as $\epsilon\rightarrow 0$ is shown in \autoref{s:discussion}.  In  \autoref{s:conclusions} we present an example showing the effect of the constants $C_t$ and $C_s$ on the convergence rate and summarize our results.

\section{Double well potential and existence of a solution}
\label{s:peridynamic model}
In this section, we present the nonlinear nonlocal model. Let $D\subset \bbR^d$, $d=2,3$ be the material domain with characteristic length-scale of unity. Let $\epsilon\in (0,1]$ be the size of horizon across which nonlocal interaction between points takes place. The material point $\bx\in D$ interacts nonlocally with all material points inside a horizon of length $\epsilon$. Let $H_{\epsilon}(\bx)$ be the ball of radius $\epsilon$ centered at $\bx$ containing all points $\by$ that interact with $\bx$. After deformation the material point $\bx$ assumes position $\bz = \bx + \bu(\bx)$. In this treatment we assume infinitesimal displacements and the strain is written in terms of the displacement $\bu$ as
\begin{align*}
S=S(\by,\bx;\bu) &:= \dfrac{\bu(\by) - \bu(\bx)}{\abs{\by - \bx}} \cdot \dfrac{\by - \bx}{\abs{\by - \bx}}.
\end{align*}

Let $W^{\epsilon}(S,\by - \bx)$ be the nonlocal potential density per unit length between material point $\by$ and $\bx$. The energy density at $\bx$ is given by
\begin{align*}
\bW^{\epsilon}(S,\bx) = \dfrac{1}{\epsilon^d \omega_d} \int_{H_{\epsilon}(\bx)} |\by-\bx|W^{\epsilon}(S, \by - \bx) d\by,
\end{align*}
where $\omega_d$ is the volume of a unit ball in $d$-dimension and $\epsilon^d \omega_d$ is the volume of the ball of radius $\epsilon$. The potential energy is written as
\begin{align*}
PD^{\epsilon}(\bu) &= \int_D \bW^{\epsilon}(S(\bu), \bx) d\bx,
\end{align*}
and the displacement field satisfies following equation of motion
\begin{align}\label{eq:per equation}
\rho \partial^2_{tt} \bu(t,\bx) &= -\del PD^{\epsilon}(\bu) + \bb(t,\bx)
\end{align}
for all $\bx \in D$. Here we have
\begin{align*}
-\del PD^{\epsilon}(\bu)(\bx) = \dfrac{2}{\epsilon^d \omega_d} \int_{H_{\epsilon}(\bx)} \partial_S W^{\epsilon}(S,\by - \bx) \dfrac{\by - \bx}{\abs{\by - \bx}} d\by
\end{align*}
where $\bb(t,\bx)$ is the body force, $\rho$ is the density and $\partial_S W^{\epsilon}$ is the derivative of potential with respect to the strain. 


We prescribe the zero Dirichlet condition on the boundary of $D$
\begin{align}\label{eq:per bc}
\bu(\bx) = \bzero \qquad \forall \bx \in \partial D,
\end{align}
where we have denoted the boundary by  $\partial D$. We extend the zero boundary condition outside $D$ to $\mathbb{R}^3$.

The peridynamic equation, boundary conditions, and  initial conditions 
\begin{align}\label{eq:per initialvalues}
\bu(0,\bx) = \bu_0(\bx) \qquad \partial_t\bu(0,\bx)=\bv_0(\bx)
\end{align}
determine the peridynamic evolution $\bu(t,\bx)$.

\paragraph{Peridynamics energy}
The total energy $\mathcal{E}^\epsilon(\bu)(t)$ is given by the sum of kinetic and potential energy given by
\begin{align}\label{eq:def energy}
\mathcal{E}^\epsilon(\bu)(t) &= \frac{1}{2} ||\dot{\bu}(t)||_{L^2(D;\bbR^d)} + PD^\epsilon(\bu(t)),
\end{align}
where potential energy $PD^\epsilon$ is given by
\begin{align*}
PD^{\epsilon}(\bu) &= \int_D \left[ \dfrac{1}{\epsilon^d \omega_d} \int_{H_\epsilon(\bx)} W^{\epsilon}(S(\bu), \by-\bx) d\by \right] d\bx.
\end{align*}
Differentiation of \autoref{eq:def energy} gives the identity
\begin{align}
\dfrac{d}{dt} \mathcal{E}^\epsilon(\bu)(t) = (\ddot{\bu}(t), \dot{\bu}(t)) - (-\del PD^\epsilon(\bu(t)), \dot{\bu}(t)) \label{eq:energy relat},
\end{align}

where $(\cdot,\cdot)$ is the inner product on $L^2(\mathbb{R}^d,D)$ and $\Vert\cdot\Vert_{L^2(\mathbb{R}^d,D)}$ is the associated norm.

\subsection{Nonlocal potential}
We consider the nonlocal two point interaction potential density $W^{\epsilon}$ of the form
\begin{align}\label{eq:per pot}
W^{\epsilon}(S, \by - \bx) &= \omega(\bx)\omega(\by)\dfrac{J^{\epsilon}(\abs{\by - \bx})}{\epsilon\abs{\by-\bx}} f(\abs{\by - \bx} S^2)
\end{align}
where $f: \bbR^{+} \to \bbR$ is assumed to be positive, smooth and concave with following properties
\begin{align}\label{eq:per asymptote}
\lim_{r\to 0^+} \dfrac{f(r)}{r} = f'(0), \qquad \lim_{r\to \infty} f(r) = f_{\infty} < \infty 
\end{align}

The potential $W^{\epsilon}(S, \by - \bx)$ is of double well type and convex near the origin where it has one well and concave and bounded at infinity where it has the second well. $J^{\epsilon}(\abs{\by - \bx})$ models the influence of separation between points $\by$ and $\bx$. We define $J^{\epsilon}$ by rescaling $J(\abs{\bxi})$, i.e. $J^{\epsilon}(\abs{\bxi}) = J(\abs{\bxi}/\epsilon)$. Here $J$ is zero outside the ball $H_1(\bzero)$ and satisfies $0\leq J(\abs{\bxi}) \leq M$ for all $\bxi \in H_1(\bzero)$. The domain function $\omega$ enforces boundary conditions on $\partial_SW^\epsilon$ at the boundary of the body $D$. Here the boundary is denoted by $\partial D$ and $\omega$ is a nonnegative differentiable function $0\leq \omega\leq 1$. On the boundary $\omega=0$ and $\omega=1$ for points $\bx$ inside $D$ with distance greater than $\epsilon$ away from the boundary. We continue $\omega$ by zero for all points outside $D$.

\begin{figure}
\centering
\includegraphics[scale=0.25]{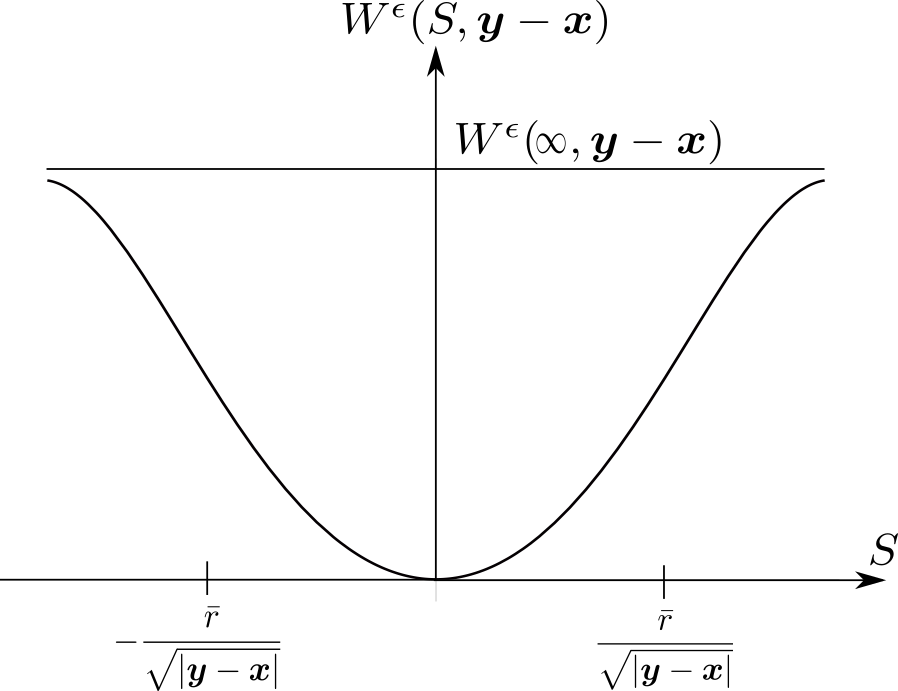}

\caption{Two point potential $W^\epsilon(S,\by - \bx)$ as a function of strain $S$ for fixed $\by - \bx$.}
 \label{fig:per pot}
\end{figure}

The potential described in \autoref{eq:per pot} gives the convex-concave dependence of $W(S,\by - \bx)$ on the strain $S$ for fixed $\by - \bx$, see \autoref{fig:per pot}. Initially the force is elastic for small strains and then softens as the strain becomes larger. The critical strain where the force between $\bx$ and $\by$ begins to soften is given by $S_c(\by, \bx) := \bar{r}/\sqrt{\abs{\by - \bx}}$ and the force decreases monotonically for
\begin{align*}
\abs{S(\by, \bx;\bu)} > S_c.
\end{align*}
Here $\bar{r}$ is the inflection point of $r:\to f(r^2)$ and is the root of following equation
\begin{align*}
f'({r}^2) + 2{r}^2 f''({r}^2) = 0.
\end{align*}

\begin{figure}
\centering
\includegraphics[scale=0.25]{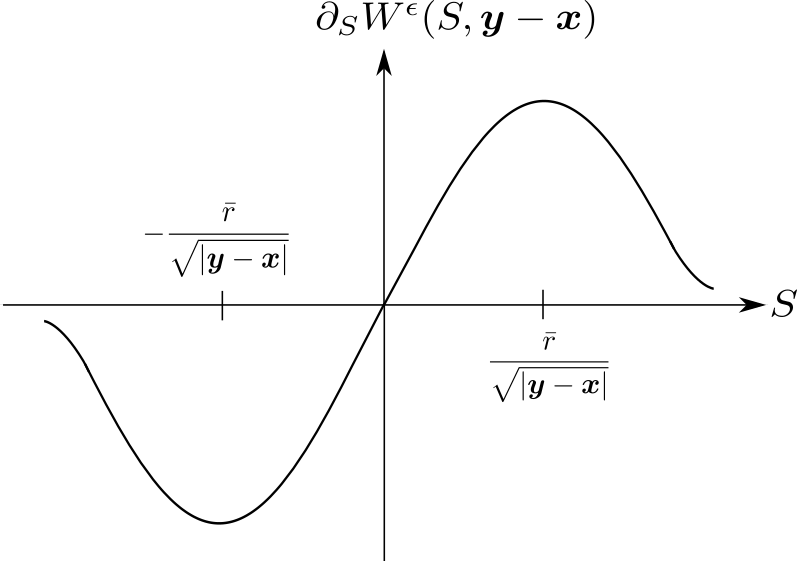}
\caption{Nonlocal force $\partial_S W^\epsilon(S,\by - \bx)$ as a function of strain $S$ for fixed $\by - \bx$. Second derivative of $W^\epsilon(S,\by-\bx)$ is zero at $\pm \bar{r}/\sqrt{|\by -\bx|}$.}
 \label{fig:first der per pot}
\end{figure}

\subsection{Existence of solution}\label{ss:existence holder}
Let $\Cholder{\perdrd}$ be the H\"{o}lder space with exponent $\gamma \in (0,1]$. 
The closure of continuous functions with compact support on $D$ in the supremum norm is denoted by $C_0(D)$. We identify functions in $C_0(D)$ with their unique continuous extensions to $\overline{D}$. It is easily seen that functions belonging to this space take the value zero on the boundary of $D$, see e.g. \cite{MA-Driver}. We introduce $C_0^{0,\gamma}(D)=C^{0,\gamma}(D)\cap C_0(D)$. In this paper we extend all functions in $C_0^{0,\gamma}(D)$ by zero outside $D$.
The norm of $\bu \in \Cholderz{\perdrd}$ is taken to be
\begin{align*}
\Choldernorm{\bu}{\perdrd} &:= \sup_{\bx \in D} \abs{\bu(\bx)} + \left[\bu \right]_{\Cholder{\perdrd}},
\end{align*}
where $\left[\bu \right]_{\Cholder{\perdrd}}$ is the H\"{o}lder semi norm and given by
\begin{align*}
\left[\bu \right]_{\Cholder{\perdrd}} &:= \sup_{\substack{\bx\neq \by,\\
\bx,\by \in D}} \dfrac{\abs{\bu(\bx)-\bu(\by)}}{\abs{\bx - \by}^\gamma},
\end{align*}
and $\Cholderz{\perdrd}$ is a Banach space with this norm. Here we make the hypothesis that the domain function $\omega$ belongs to $\Cholderz{\perdrd}$.

We write the evolution \autoref{eq:per equation} as an equivalent first order system with $y_1(t)=\bu(t)$ and $y_2(t)=\bv(t)$ with $\bv(t)=\partial_t\bu(t)$. Let $y = (y_1, y_2)^T$ where $y_1,y_2 \in \Cholderz{\perdrd}$ and let $F^{\epsilon}(y,t) = (F^{\epsilon}_1(y,t), F^{\epsilon}_2(y,t))^T$ such that
\begin{align}
F^\epsilon_1(y,t) &:= y_2 \label{eq:per first order eqn 1} \\
F^\epsilon_2(y, t) &:= -\del PD^{\epsilon}(y_1) + \bb(t). \label{eq:per first order eqn 2}
\end{align}
The initial boundary value associated with the evolution \autoref{eq:per equation} is equivalent to the initial boundary  value problem for the first order system given by
\begin{align}\label{eq:per first order}
\dfrac{d}{dt}y = F^{\epsilon}(y,t),
\end{align}
with initial condition given by $y(0) = (\bu_0, \bv_0)^T \in \Cholderz{\perdrd}\times\Cholderz{\perdrd}$.

The function $F^{\epsilon}(y,t)$ satisfies the Lipschitz continuity given by the following theorem.

{\vskip 2mm}
\begin{proposition}\label{prop:lipschitz}
\textbf{Lipschitz continuity and bound}\\
Let $X = \Cholderz{\perdrd} \times \Cholderz{\perdrd}$. The function $F^\epsilon(y,t) = (F^\epsilon_1, F^\epsilon_2)^T$, as defined in \autoref{eq:per first order eqn 1} and \autoref{eq:per first order eqn 2}, is Lipschitz continuous in any bounded subset of $X$. We have, for any $y,z \in X$ and $t> 0$,
\begin{align}\label{eq:lipschitz property of F}
&\normX{F^{\epsilon}(y,t) - F^{\epsilon}(z,t)}{X} \notag \\
&\leq \dfrac{\left( L_1 + L_2 \left( \Vert \omega\Vert_{C^{0,\gamma}(D)}+\normX{y}{X} + \normX{z}{X} \right)  \right)}{\epsilon^{2 + \alpha(\gamma)}} \normX{y-z}{X}
\end{align}
where $L_1, L_2$ are independent of $\bu,\bv$ and depend on peridynamic potential function $f$ and influence function $J$ and the exponent $\alpha(\gamma)$ is given by
\begin{align*}
\alpha(\gamma) = \begin{cases}
0 &\qquad \text{if }\gamma \geq 1/2 \\
1/2 - \gamma &\qquad \text{if } \gamma < 1/2 .
\end{cases}
\end{align*}
Furthermore for any $y \in X$ and any $t\in [0,T]$, we have the bound
\begin{align}\label{eq:bound on F}
\normX{F^\epsilon(y,t)}{X} &\leq \dfrac{L_3}{\epsilon^{2+\alpha(\gamma)}} (1+\Vert \omega\Vert_{C^{0,\gamma}(D)} + \normX{y}{X}) + b
\end{align}
where $b = \sup_{t} \Choldernorm{\bb(t)}{\perdrd}$ and $L_3$ is independent of $y$. 
\end{proposition}
{\vskip 2mm}

\setcounter{theorem}{1}
We easily see that on choosing $z=0$ in \autoref{eq:lipschitz property of F} that $-\nabla PD^\epsilon(\bu)(\bx)$ is in $C^{0,\gamma}(D;\mathbb{R}^3)$ provided that $\bu$ belongs to $C^{0,\gamma}(D;\mathbb{R}^3)$. Since  $-\nabla PD^\epsilon(\bu)(\bx)$ takes the value $0$ on $\partial D$ we conclude that $-\nabla PD^\epsilon(\bu)(\bx)$ belongs to  $C^{0,\gamma}_0(D;\mathbb{R}^3)$.

In Theorem 6.1 of \cite{CMPer-Lipton}, the Lipschitz property of a peridynamic force is shown in $X = \Ltwo{\perdrd} \times \Ltwo{\perdrd}$. It is given by
\begin{align}\label{eq: lipshitz}
\normX{F^{\epsilon}(y,t) - F^{\epsilon}(z,t)}{X} &\leq \dfrac{L}{\epsilon^2}\normX{y - z}{X} \qquad \forall y,z \in X, \forall t\in [0,T]
\end{align}
for all $y,z \in \Ltwoz{\perdrd}^2$. For this case $L$ does not depend on $\bu, \bv$. We now state the existence theorem.

The following theorem gives the existence and uniqueness of solution in any given time domain $I_0 = (-T, T)$.

{\vskip 2mm}
\begin{theorem}\label{thm:existence over finite time domain} 
\textbf{Existence and uniqueness of H\"older solutions of cohesive dynamics over finite time intervals}\\
For any initial condition $x_0\in X ={ \Cholderz{\perdrd} \times \Cholderz{\perdrd}}$, time interval $I_0=(-T,T)$, and right hand side $\bb(t)$  continuous in time for $t\in I_0$ such that $\bb(t)$ satisfies $\sup_{t\in I_0} {||\bb(t)||_{\Cholder{}}}<\infty$, there is a unique solution $y(t)\in C^1(I_0;X)$ of
\begin{equation*}
y(t)=x_0+\int_0^tF^\epsilon(y(\tau),\tau)\,d\tau,
\label{10}
\end{equation*}
or equivalently
\begin{equation*}
y'(t)=F^\epsilon(y(t),t),\hbox{with    $y(0)=x_0$},
\label{11}
\end{equation*}
where $y(t)$ and $y'(t)$ are Lipschitz continuous in time for $t\in I_0$.
\end{theorem}
{\vskip 2mm}

The proof of this theorem is given in \autoref{s:proof existence}. We now describe the finite difference scheme and analyze its convergence to H\"older continuous solutions of cohesive dynamics.

\section{Finite difference approximation}
\label{s:finite difference}
In this section, we present the finite difference scheme and compute the rate of convergence. We first consider the semi-discrete approximation and prove the bound on energy of semi-discrete evolution in terms of initial energy and the work done by body forces. 

Let $h$ be the size of a mesh and $\Delta t$ be the size of time step. We will keep $\epsilon$ fixed and assume that $h< \epsilon<1$. Let $D_h = D\cap(h \bbZ)^d$ be the discretization of material domain. Let $i\in \bbZ^d$ be the index such that $\bx_i = hi \in D$. Let $U_i$ is the unit cell of volume $h^d$ corresponding to the grid point $\bx_i$. The exact solution evaluated at grid points is denoted by $(\bu_i(t),\bv_i(t))$.

\subsection{Time discretization}
Let $[0,T] \cap (\Delta t \bbZ)$ be the discretization of time domain where $\Delta t$ is the size of time step. Denote fully discrete solution at $(t^k = k\Delta t, \bx_i = ih)$ as $(\hat{\bu}^k_{i}, \hat{\bv}^k_i)$. Similarly, the exact solution evaluated at grid points is denoted by $(\bu^k_i,\bv^k_i)$. We enforce boundary condition $\hat{\bu}^k_i = \bzero$ for all $\bx_i \notin D$ and for all $k$.

We begin with the forward Euler time discretization, with respect to velocity, and the finite difference scheme for $(\hat{\bu}^k_{i}, \hat{\bv}^k_i)$ is written
\begin{align}
\dfrac{\buhat^{k+1}_i - \buhat^k_i}{\Delta t} &= \bvhat^{k+1}_i \label{eq:finite diff eqn u} \\
\dfrac{\bvhat^{k+1}_i - \bvhat^k_i}{\Delta t} &= -\del PD^{\epsilon}(\buhat^k)(\bx_i) + \bb^k_i \label{eq:finite diff eqn v}
\end{align}
The scheme is complemented with the discretized initial conditions $\buhat^{0}_i =(\buhat_0)_i$ and $\bvhat^{0}_i =(\bvhat_0)_i$. If we substitute \autoref{eq:finite diff eqn u} into \autoref{eq:finite diff eqn v}, we get standard Central difference scheme in time for second order in time differential equation. Here we have assumed, without loss of generality, $\rho = 1$.

The piecewise constant extensions of the discrete sets $\{\hat{\bu}^k_i\}_{i\in \bbZ^d}$ and $\{\hat{\bv}^k_i\}_{i\in \bbZ^d}$ are given by
\begin{align*}
\hat{\bu}^k(\bx) &:= \sum_{i, \bx_i \in D} \hat{\bu}^k_i \chi_{U_i}(\bx) \\
\hat{\bv}^k(\bx) &:= \sum_{i, \bx_i \in D} \hat{\bv}^k_i \chi_{U_i}(\bx)
\end{align*}
In this way we represent the finite difference solution as a piecewise constant function. We will show this function provides an $L^2$ approximation of the exact solution.

\begin{figure}[h]
\centering
\includegraphics[scale=0.5]{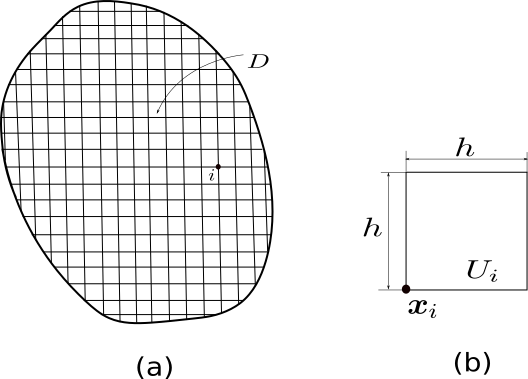}
\caption{(a) Typical mesh of size $h$. (b) Unit cell $U_i$ corresponding to material point $\bx_i$.}\label{fig:peridynamic mesh}
\end{figure}

\subsubsection{Convergence results}
In this section we provide upper bounds on the rate of convergence of the discrete approximation to the solution of the peridynamic evolution. The $L^2$ approximation error $E^k$ at time $t^k$,  for $0<t^k\leq T$ is defined as 
\begin{align*}
E^k &:= \Ltwonorm{\buhat^k - \bu^k}{\perdrd} + \Ltwonorm{\bvhat^k- \bv^k}{\perdrd} 
\end{align*}
The upper bound on the convergence rate of the approximation error is given by the following theorem.

{\vskip 2mm}
\begin{theorem}\label{thm:convergence}
\textbf{Convergence of finite difference approximation (forward Euler time discretization)}\\
Let $\epsilon>0$ be fixed. Let $(\bu, \bv)$ be the solution of peridynamic equation \autoref{eq:per first order}. We assume $\bu, \bv \in \Ctwointime{[0,T]}{\Cholderz{\perdrd}}$.   Then the finite difference scheme given by \autoref{eq:finite diff eqn u} and \autoref{eq:finite diff eqn v} is consistent in both  time and spatial discretization and converges to the exact solution uniformly in time with respect to the $\Ltwo{\perdrd}$ norm. If we assume the error at the initial step is zero then the error $E^k$ at time $t^k$ is bounded and to leading order in the time step $\Delta t$ satisfies
\begin{align}\label{eq: first est}
\sup_{0\leq k \leq T/\Delta t} E^k\leq O\left( C_t\Delta t + C_s\dfrac{h^\gamma}{\epsilon^2} \right),
\end{align}
where constant $C_s$ and $C_t$ are independent of $h$ and $\Delta t$ and $C_s$ depends on the H\"older norm of the solution and $C_t$ depends on the $L^2$ norms of time derivatives of the solution. 
\end{theorem}
{\vskip 2mm}
Here we have assumed the initial error to be zero for ease of exposition only. 

We remark that the explicit constants leading to \autoref{eq: first est} can be large. The inequality that delivers \autoref{eq: first est}  is given to leading order by
\begin{align}\label{eq: fund est initial}
\sup_{0\leq k \leq T/\Delta t} E^k\leq \exp \left[T (1 + 6\bar{C}/\epsilon^2) \right] T \left[ C_t \Delta t + (C_s/\epsilon^2) h^\gamma \right],
\end{align}
where the constants $\bar{C}$, $C_t$ and $C_s$ are given by \autoref{eq: bar C},  \autoref{eq:const Ct}, and \autoref{eq:const Cs}.  The explicit constant $C_t$ depends on the spatial $L^2$ norm of the time derivatives of the solution and $C_s$ depends on the spatial H\"older continuity of the solution  and  the constant $\bar{C}$. This constant is bounded independently of horizon $\epsilon$.  Although the constants are necessarily pessimistic they deliver a-priori error estimates  and an example is discussed in  \autoref{s:conclusions}. 

An identical convergence rate can be established for the general one step scheme and we state it below.

{\vskip 2mm}
\begin{theorem}\label{thm:convergence general}
\textbf{Convergence of finite difference approximation (General single step time discretization)}\\
Let us assume that the hypothesis of \autoref{thm:convergence} holds. Fix $\theta \in [0,1]$, and let $(\buhat^k , \bvhat^k)^T$ be the solution of following finite difference equation
\begin{align}
\dfrac{\buhat^{k+1}_i - \buhat^k_i}{\Delta t} &= (1 - \theta) \bvhat^k_i + \theta \bvhat^{k+1}_i \label{eq:finite diff eqn u general} \\
\dfrac{\bvhat^{k+1}_i - \bvhat^k_i}{\Delta t} &= (1-\theta)\left( -\del PD^{\epsilon}(\buhat^k)(\bx_i) + \bb^k_i\right) + \theta \left( - \del PD^{\epsilon}(\buhat^{k+1})(\bx_i) + \bb^{k+1}_i \right).  \label{eq:finite diff eqn v general}
\end{align}
Then, for any fixed $\theta\in [0,1]$, there exists a constant $K>0$ independent of $(\buhat^k , \bvhat^k)^T$ and $(\bu^k , \bv^k)^T$, such that for $\Delta t<K{\epsilon^2}$ the finite difference scheme given by \autoref{eq:finite diff eqn u general} and \autoref{eq:finite diff eqn v general} is consistent in time and spatial discretization. if we assume the error at the initial step is zero then the error $E^k$ at time $t^k$ is bounded and satisfies
\begin{align*}
\sup_{0\leq k \leq T/\Delta t} E^k \leq O\left( C_t\Delta t + C_s\dfrac{h^\gamma}{\epsilon^2} \right).
\end{align*}

The constant $K$ is given by the explicit formula $K=1/\bar{C}$ where $\bar{C}$ is described by equation \autoref{eq: bar C}. 
Furthermore for the Crank Nicholson scheme, $\theta = 1/2$, if we assume the solutions $\bu,\bv$ belong to  $\Cthreeintime{[0,T]}{\Cholderz{\perdrd}}$, then the approximation error $E^k$ satisfies
\begin{align*}
\sup_{0\leq k \leq T/\Delta t} E^k \leq O\left(\bar{C}_t (\Delta t)^2 + C_s\dfrac{h^\gamma}{\epsilon^2} \right),
\end{align*}
where $\bar{C}_t$ is independent of $\Delta t$ and $h$ and is given by \autoref{eq: timebound}.
\end{theorem}
{\vskip 2mm}
As before we assume that the error in the initial data is zero for ease of exposition. The proofs of \autoref{thm:convergence} and \autoref{thm:convergence general} are given in the following sections.

\textbf{Remark.} {\em In \autoref{thm:convergence general}, we have stated a condition on $\Delta t$ for which the convergence estimate holds. This condition naturally occurs in the analysis and is related to the Lipschitz continuity of the peridnamic force with respect to the $L^2$ norm, see \autoref{eq: lipshitzl2}}.

\subsubsection{Error analysis}\label{ss:error analysis1}

\autoref{thm:convergence} and \autoref{thm:convergence general} are proved along similar lines. In both cases we define the $L^2$-projections of the actual solutions onto the space of piecewise constant functions defined over the cells $U_i$. These are given as follows. Let $(\butilde^k_i, \bvtilde^k_i)$ be the average of the exact solution $(\bu^k, \bv^k)$ in the unit cell $U_i$ given by
\begin{align*}
\butilde^k_i &:= \dfrac{1}{h^d} \int_{U_i} \bu^k(\bx) d\bx \\
\bvtilde^k_i &:= \dfrac{1}{h^d} \int_{U_i} \bv^k(\bx) d\bx 
\end{align*}
and the $L^2$ projection of the solution onto piecewise  constant functions are  $(\butilde^k, \bvtilde^k)$ given by
\begin{align}
\butilde^k(\bx) &:= \sum_{i, \bx_i \in D} \butilde^k_i \chi_{U_i}(\bx) \label{eq:periodpiecewise ext1} \\
\bvtilde^k(\bx) &:= \sum_{i, \bx_i \in D}  \bvtilde^k_i\chi_{U_i}(\bx) \label{eq:periodpiecewise ext2}
\end{align}

The error between $(\buhat^k, \bvhat^k)^T$ with $(\bu(t^k), \bv(t^k))^T$ is now split into two parts. From the triangle inequality, we have

\begin{align*}
\Ltwonorm{\buhat^k - \bu(t^k)}{\perdrd} &\leq \Ltwonorm{\buhat^k - \butilde^k}{\perdrd} + \Ltwonorm{\butilde^k - \bu^k}{\perdrd}  \\
\Ltwonorm{\bvhat^k - \bv(t^k)}{\perdrd} &\leq \Ltwonorm{\bvhat^k - \bvtilde^k}{\perdrd} + \Ltwonorm{\bvtilde^k - \bv^k}{\perdrd} 
\end{align*}

In \autoref{ss:error analysis} and \autoref{ss:implicit} we will show that the error between the $L^2$ projections of the actual solution and the discrete approximation for both forward Euler and implicit one step methods decay according to
\begin{align}\label{eq:write estimate error ek}
\sup_{0\leq k \leq T/\Delta t} \left( \Ltwonorm{\buhat^k - \butilde^k}{\perdrd} + \Ltwonorm{\bvhat^k - \bvtilde^k}{\perdrd} \right)  &= O\left( \Delta t + \dfrac{h^\gamma}{\epsilon^2} \right).
\end{align}
In what follows we can estimate the terms
\begin{align}\label{eq: per convgest}
\Ltwonorm{\butilde^k - \bu(t^k)}{}\hbox{  and    }
\Ltwonorm{\bvtilde^k - \bv(t^k)}{}
\end{align}
and show they  go to zero at a rate of $h^\gamma$ uniformly in time. The estimates given by \autoref{eq:write estimate error ek} together with the $O(h^\gamma)$ estimates for \autoref{eq: per convgest} establish \autoref{thm:convergence} and \autoref{thm:convergence general}. We now establish the $L^2$ estimates for the differences $\butilde^k - \bu(t^k)$ and $\bvtilde^k - \bv(t^k)$. 

We  write
\begin{align}
&\Ltwonorm{\butilde^k - \bu^k}{\perdrd}^2  \notag \\
&= \sum_{i, \bx_i \in D} \int_{U_i} \abs{\butilde^k(\bx) - \bu^k(\bx)}^2 d\bx \notag \\
&= \sum_{i,\bx_i \in D} \int_{U_i} \abs{ \dfrac{1}{h^d} \int_{U_i} (\bu^k(\by) - \bu^k( \bx)) d\by }^2 d\bx \notag \\
&= \sum_{i,\bx_i \in D} \int_{U_i} \left[ \dfrac{1}{h^{2d}} \int_{U_i} \int_{U_i} (\bu^k(\by) - \bu^k( \bx)) \cdot (\bu^k(\bz) - \bu^k( \bx)) d\by d\bz \right] d\bx \notag \\
&\leq \sum_{i,\bx_i \in D} \int_{U_i} \left[ \dfrac{1}{h^d} \int_{U_i} \abs{\bu^k(\by) - \bu^k(\bx)}^2 d\by \right] d\bx \label{eq:estimate butilde and bu}
\end{align}
where we used Cauchy's inequality and Jensen's inequality. For $\bx,\by \in U_i$, $\abs{\bx - \by} \leq c h$, where $c = \sqrt{2}$ for $d=2$ and $c=\sqrt{3}$ for $d=3$. Since $\bu \in \Cholderz{}$ we have
\begin{align}\label{eq:estimate bu holder}
\abs{\bu^k(\bx) - \bu^k(\by)} &= \abs{\bx - \by}^\gamma \dfrac{\abs{\bu^k(\by) -\bu^k(\bx)}}{\abs{\bx - \by}^\gamma}  \notag \\
&\leq c^{\gamma} h^\gamma \Choldernorm{\bu^k}{\perdrd} \leq c^{\gamma} h^\gamma \sup_t \Choldernorm{\bu(t)}{\perdrd}
\end{align}
and substitution in \autoref{eq:estimate butilde and bu} gives
\begin{align*}
\Ltwonorm{\butilde^k - \bu^k}{\perdrd}^2  &\leq c^{2\gamma} h^{2\gamma} \sum_{i,\bx_i \in D} \int_{U_i} d\bx \left( \sup_t \Choldernorm{\bu(t)}{\perdrd} \right)^2 \notag \\
&\leq c^{2\gamma}  |D|  h^{2\gamma} \left( \sup_t \Choldernorm{\bu(t)}{\perdrd} \right)^2 .
\end{align*}

A similar estimate can be derived for $||\bvtilde^k - \bv^k||_{L^2}$ and substitution of the  estimates into \autoref{eq: per convgest} gives
\begin{align*}
\sup_{k} \left( \Ltwonorm{\butilde^k - \bu(t^k)}{\perdrd}  + \Ltwonorm{\bvtilde^k - \bv(t^k)}{\perdrd} \right) = O(h^\gamma).
\end{align*}

In the next section we establish the error estimate \autoref{eq:write estimate error ek} for both forward Euler and general one step schemes in \autoref{ss:error analysis} and \autoref{ss:implicit}.

\subsubsection{Error analysis for approximation of $L^2$ projection of the the exact solution}\label{ss:error analysis}
In this sub-section, we estimate the difference between approximate solution $(\buhat^k,\bvhat^k)$ and the $L^2$ projection of the exact solution onto piece wise constant functions given by $(\butilde^k, \bvtilde^k)$, see \autoref{eq:periodpiecewise ext1} and \autoref{eq:periodpiecewise ext2} . Let the differences be denoted by $\be^k(u) := \buhat^k - \butilde^k$ and  $\be^k(v ):= \bvhat^k - \bvtilde^k$ and their evaluation at grid points are $\be^k_i(u) := \buhat^k_i - \butilde^k_i$ and $\be^k_i(v) := \bvhat^k_i - \bvtilde^k_i$. 
Subtracting $(\butilde^{k+1}_i - \butilde^k_i)/\Delta t$ from \autoref{eq:finite diff eqn u} gives
\begin{align*}
& \dfrac{\buhat^{k+1}_i - \buhat^k_i}{\Delta t} - \dfrac{\butilde^{k+1}_i - \butilde^k_i}{\Delta t}  \\
&= \bvhat^{k+1}_i - \dfrac{\butilde^{k+1}_i - \butilde^k_i}{\Delta t} \notag \\
&= \bvhat^{k+1}_i - \bvtilde^{k+1}_i + \left( \bvtilde^{k+1}_i - \dparder{\butilde^{k+1}_i}{t} \right) + \left( \dparder{\butilde^{k+1}_i}{t} - \dfrac{\butilde^{k+1}_i - \butilde^k_i}{\Delta t} \right).
\end{align*}
Taking the average over unit cell $U_i$ of the exact peridynamic equation \autoref{eq:per first order} at time $t^k$, we will get $\bvtilde^{k+1}_i - \dparder{\butilde^{k+1}_i}{t}  = 0$. Therefore, the equation for $\be^k_i(u)$ is given by
\begin{align}\label{eq:error eqn in u}
\be^{k+1}_i(u) = \be^k_i(u) + \Delta t \be^{k+1}_i(v) + \Delta t\tau^{k}_i(u),
\end{align}
where we identify the discretization error as
\begin{align}\label{eq:consistency error in u}
\tau^k_i(u) &:= \dparder{\butilde^{k+1}_i}{t} - \dfrac{\butilde^{k+1}_i - \butilde^k_i}{\Delta t}.
\end{align}

Similarly, we subtract $(\bvtilde^{k+1}_i - \bvtilde^k_i)/\Delta t$ from \autoref{eq:finite diff eqn v} and add and subtract terms to get
\begin{align}\label{eq:error in v 1}
\dfrac{\bvhat^{k+1}_i - \bvhat^k_i}{\Delta t} - \dfrac{\bvtilde^{k+1}_i - \bvtilde^k_i}{\Delta t} &= - \del PD^{\epsilon}(\buhat^k)(\bx_i) + \bb^k_i - \dparder{\bv^k_i}{t} + \left( \dparder{\bv^k_i}{t} - \dfrac{\bvtilde^{k+1}_i - \bvtilde^k_i}{\Delta t}\right) \notag \\
&= - \del PD^{\epsilon}(\buhat^k)(\bx_i) + \bb^k_i - \dparder{\bv^k_i}{t} \notag \\
&\quad + \left( \dparder{\bvtilde^k_i}{t} - \dfrac{\bvtilde^{k+1}_i - \bvtilde^k_i}{\Delta t}\right) + \left( \dparder{\bv^k_i}{t} - \dparder{\bvtilde^k_i}{t}\right),
\end{align}
where we identify $\tau^k_i(v)$ as follows
\begin{align}\label{eq:consistency error in v}
\tau^k_i(v) &:= \dparder{\bvtilde^k_i}{t} - \dfrac{\bvtilde^{k+1}_i - \bvtilde^k_i}{\Delta t}.
\end{align}
Note that in $\tau^k(u)$ we have $\dparder{\butilde^{k+1}_i}{t}$ and from the exact peridynamic equation, we have
\begin{align}\label{eq:exact per eqn v 1}
\bb^k_i - \dparder{\bv^k_i}{t} = \del PD^{\epsilon}(\bu^k)(\bx_i).
\end{align}
Combining \autoref{eq:error in v 1}, \autoref{eq:consistency error in v}, and \autoref{eq:exact per eqn v 1}, to get
\begin{align*}
\be^{k+1}_i(v) &= \be^k_i(v) + \Delta t \tau^k_i(v) + \Delta t \left( \dparder{\bv^k_i}{t} - \dparder{\bvtilde^k_i}{t}\right)  \notag \\
&\quad + \Delta t \left( -\del PD^{\epsilon}(\buhat^k)(\bx_i) + \del PD^{\epsilon}(\bu^k)(\bx_i) \right) \notag \\
&= \be^k_i(v) + \Delta t \tau^k_i(v) + \Delta t \left( \dparder{\bv^k_i}{t} - \dparder{\bvtilde^k_i}{t}\right) \notag \\
&\quad + \Delta t \left( -\del PD^{\epsilon}(\buhat^k)(\bx_i) + \del PD^{\epsilon}(\butilde^k)(\bx_i) \right) \notag \\
&\quad + \Delta t \left( -\del PD^{\epsilon}(\butilde^k)(\bx_i) + \del PD^{\epsilon}(\bu^k)(\bx_i) \right).
\end{align*}
The spatial discretization error $\sigma^k_i(u)$ and $\sigma^k_i(v)$ is given by
\begin{align}
\sigma^k_i(u) &:= \left( -\del PD^{\epsilon}(\butilde^k)(\bx_i) + \del PD^{\epsilon}(\bu^k)(\bx_i)  \right) \label{eq:consistency error in u spatial} \\
\sigma^k_i(v) &:=  \dparder{\bv^k_i}{t} - \dparder{\bvtilde^k_i}{t}. \label{eq:consistency error in v spatial}
\end{align}
We finally have 
\begin{align}\label{eq:error eqn in v}
\be^{k+1}_i(v) &= \be^k_i(v) + \Delta t \left(\tau^k_i(v) + \sigma^k_i(u) + \sigma^k_i(v) \right) \notag \\
&\quad + \Delta t \left( -\del PD^{\epsilon}(\buhat^k)(\bx_i) + \del PD^{\epsilon}(\butilde^k)(\bx_i) \right).
\end{align}
We now show the consistency and stability properties of the numerical scheme.

\subsubsection{Consistency}\label{sss:consistency}
We deal with the error in time discretization and the error in spatial discretization error separately. The time discretization error follows easily using the Taylor's series while spatial the discretization error uses properties of the nonlinear peridynamic force.

\

\textbf{Time discretization: }We first estimate the time discretization error. A Taylor series expansion is used to estimate $\tau^k_i(u)$ as follows
\begin{align*}
\tau^k_i(u) &= \dfrac{1}{h^d} \int_{U_i} \left( \dparder{\bu^k(\bx)}{t} - \dfrac{\bu^{k+1}(\bx) - \bu^k(\bx)}{\Delta t} \right) d\bx \\
&= \dfrac{1}{h^d} \int_{U_i} \left( -\dfrac{1}{2} \dsecder{\bu^k(\bx)}{t} \Delta t + O((\Delta t)^2) \right) d\bx .
\end{align*}
Computing the $\Ltwo{}$ norm of $\tau^k_i(u)$ and using Jensen's inequality gives
\begin{align*}
\Ltwonorm{\tau^k(u)}{\perdrd} &\leq \frac{\Delta t}{2} \Ltwonorm{\dsecder{\bu^k}{t}}{\perdrd} + O((\Delta t)^2) \notag \\
&\leq \frac{\Delta t}{2} \sup_{t} \Ltwonorm{\dsecder{\bu(t)}{t}}{\perdrd} + O((\Delta t)^2).
\end{align*}
Similarly, we have
\begin{align*}
\Ltwonorm{\tau^k(v)}{\perdrd} = \frac{\Delta t}{2} \sup_{t} \Ltwonorm{\dsecder{\bv(t)}{t}}{\perdrd} + O((\Delta t)^2).
\end{align*}

\

\textbf{Spatial discretization: }We now estimate the spatial discretization error. Substituting the definition of $\bvtilde^k$ and following the similar steps employed in \autoref{eq:estimate bu holder}, gives
\begin{align*}
\abs{\sigma^k_i(v)} &= \abs{\dparder{\bv^k_i}{t} - \dfrac{1}{h^d}\int_{U_i} \dparder{\bv^k(\bx)}{t} d\bx} \leq c^\gamma h^{\gamma} \int_{U_i} \dfrac{1}{\abs{\bx_i - \bx}^\gamma} \abs{\dparder{\bv^k(\bx_i)}{t} - \dparder{\bv^k(\bx)}{t}} d\bx \notag \\
&\leq c^\gamma h^{\gamma} \Choldernorm{\dparder{\bv^k}{t}}{\perdrd} \leq c^\gamma h^{\gamma} \sup_{t} \Choldernorm{\dparder{\bv(t)}{t}}{\perdrd}.
\end{align*}
Taking the $\Ltwo{}$ norm of error $\sigma^k_i(v)$ and substituting the estimate above delivers
\begin{align*}
\Ltwonorm{\sigma^k(v)}{\perdrd} &\leq h^{\gamma} c^{\gamma} \sqrt{ \abs{D} } \sup_{t} \Choldernorm{\dparder{\bv(t)}{t}}{\perdrd}.
\end{align*}

Now we estimate $\abs{\sigma^k_i(u)}$. We use the notation $\bar{\bu}^k(\bx):= \bu^k(\bx+\epsilon \bxi) - \bu^k(\bx)$ and $\overline{\tilde{\bu}}^k(\bx):= \tilde{\bu}(\bx+\epsilon \bxi) -\tilde {\bu}^k(\bx)$ and choose $\bu=\bu^k$ and $\bv=\tilde{\bu}^k$ in \autoref{eq:diff force bound 1} to find that
\begin{align}\label{eq:estimate sigma u 1}
\abs{\sigma^k_i(u)} &= \abs{-\del PD^{\epsilon}(\butilde^k)(\bx_i) + \del PD^{\epsilon}(\bu^k)(\bx_i)} \notag \\
&\leq \dfrac{2C_2}{\epsilon \omega_d} \abs{ \int_{H_1(\bzero)} J(\abs{\bxi}) \dfrac{\abs{\bu^k(\bx_i + \epsilon \bxi) - \butilde^k(\bx_i + \epsilon \bxi) - (\bu^k(\bx_i) - \butilde^k(\bx_i))}}{\epsilon \abs{\bxi}} d\bxi }.
\end{align}
Here $C_2$ is the maximum of the second derivative of the profile describing the potential given by \autoref{C2}.
Following the earlier analysis, see \autoref{eq:estimate bu holder}, we find that
\begin{align}
\abs{\bu^k(\bx_i + \epsilon \bxi) - \butilde^k(\bx_i + \epsilon \bxi)} &\leq c^\gamma h^\gamma \sup_t \Choldernorm{\bu(t)}{\perdrd} \notag \\
\abs{ \bu^k(\bx_i) - \butilde^k(\bx_i) } &\leq c^\gamma h^\gamma \sup_t \Choldernorm{\bu(t)}{\perdrd}. \notag
\end{align}
For reference, we define the constant
\begin{align}\label{eq: bar C}
&\bar{C}=\frac{C_2}{\omega_d}\int_{H_1(\bzero)}J(\abs{\bxi})\dfrac{1}{\abs{\bxi}}\,d\bxi.
\end{align}

We now focus on \autoref{eq:estimate sigma u 1}. We substitute the above two inequalities to get
\begin{align*}
\abs{\sigma^k_i(u)} &\leq \dfrac{2C_2}{\epsilon^2 \omega_d} \vline \int_{H_1(\bzero)} J(\abs{\bxi}) \dfrac{1}{\abs{\bxi}} \notag \\
&\qquad \left(\abs{\bu^k(\bx_i + \epsilon \bxi) - \butilde^k(\bx_i + \epsilon \bxi)} + \abs{ \bu^k(\bx_i) - \butilde^k(\bx_i) }  \right) d\bxi \vline \notag \\
&\leq 4 h^\gamma c^\gamma \frac{\bar{C}}{\epsilon^2}  \sup_t  \Choldernorm{\bu(t)}{\perdrd}.
\end{align*}

Therefore, we have
\begin{align*}
\Ltwonorm{\sigma^k(u)}{\perdrd} &\leq h^\gamma \left( 4 c^\gamma \sqrt{|D|} \frac{\bar{C}}{\epsilon^2}  \sup_t  \Choldernorm{\bu(t)}{\perdrd} \right).
\end{align*}
This completes the proof of consistency of numerical approximation.

\subsubsection{Stability}\label{sss:stability}
Let $e^k$ be the total error at the $k^{\text{th}}$ time step. It is defined as
\begin{align*}
e^k &:=  \Ltwonorm{\be^k(u)}{\perdrd} + \Ltwonorm{\be^k(v)}{\perdrd}.
\end{align*}
To simplify the calculations, we define new term $\tau$ as 
\begin{align*}
\tau &:= \sup_t \left(\Ltwonorm{\tau^k(u)}{\perdrd}  + \Ltwonorm{\tau^k(v)}{\perdrd} \right. \notag \\
&\quad \left. + \Ltwonorm{\sigma^k(u)}{\perdrd} + \Ltwonorm{\sigma^k(v)}{\perdrd}\right).
\end{align*}
From our consistency analysis, we know that to leading order
\begin{align}\label{eq:estimate tau}
\tau &\leq C_t \Delta t + \dfrac{C_s}{\epsilon^2} h^\gamma
\end{align}
where,
\begin{align}
C_t &:= \frac{1}{2} \sup_{t} \Ltwonorm{\dsecder{\bu(t)}{t}}{\perdrd} + \frac{1}{2} \sup_{t} \Ltwonorm{\dfrac{\partial^3 \bu(t)}{\partial t^3}}{\perdrd}, \label{eq:const Ct} \\
C_s &:= c^\gamma \sqrt{|D|} \left[ \epsilon^2 \sup_{t} \Choldernorm{\dfrac{\partial^2 \bu(t)}{\partial t^2}}{\perdrd} + 4 \bar{C} \sup_t  \Choldernorm{\bu(t)}{\perdrd} \right]. \label{eq:const Cs}
\end{align}

We take $\Ltwo{}$ norm of \autoref{eq:error eqn in u} and \autoref{eq:error eqn in v} and add them. Noting the definition of $\tau$ as above, we get
\begin{align}\label{eq:error k ineq 1}
e^{k+1} &\leq e^k + \Delta t \Ltwonorm{\be^{k+1}(v)}{\perdrd} + \Delta t \tau \notag \\
&\quad + \Delta t \left( \sum_{i} h^d \abs{ -\del PD^{\epsilon}(\buhat^k)(\bx_i) + \del PD^{\epsilon}(\butilde^k)(\bx_i)}^2 \right)^{1/2}.
\end{align}

We only need to estimate the last term in above equation. Similar to the \autoref{eq:estimate sigma u 1}, we have
\begin{align*}
&\abs{ -\del PD^{\epsilon}(\buhat^k)(\bx_i) + \del PD^{\epsilon}(\butilde^k)(\bx_i)} \notag \\
&\leq \dfrac{2C_2}{\epsilon^2 \omega_d} \vline \int_{H_1(\bzero)} J(\abs{\bxi}) \dfrac{1}{\abs{\bxi}} \abs{\buhat^k(\bx_i + \epsilon \bxi) - \butilde^k(\bx_i + \epsilon \bxi) - (\buhat^k(\bx_i) - \butilde^k(\bx_i)) } d\bxi \vline \notag \\
&= \dfrac{2C_2}{\epsilon^2 \omega_d} \vline \int_{H_1(\bzero)} J(\abs{\bxi}) \dfrac{1}{\abs{\bxi}} \abs{\be^k(u)(\bx_i + \epsilon \bxi) - \be^k(u)(\bx_i)} d\bxi \vline  \notag \\
&\leq \dfrac{2C_2}{\epsilon^2 \omega_d} \vline \int_{H_1(\bzero)} J(\abs{\bxi}) \dfrac{1}{\abs{\bxi}} \left(\abs{\be^k(u)(\bx_i + \epsilon \bxi)} + \abs{\be^k(u)(\bx_i )} \right) d\bxi \vline .
\end{align*}
By $\be^k(u)(\bx)$ we mean evaluation of piecewise extension of set $\{ \be^k_i(u) \}_i$ at $\bx$. We proceed further as follows
\begin{align*}
&\abs{ -\del PD^{\epsilon}(\buhat^k)(\bx_i) + \del PD^{\epsilon}(\butilde^k)(\bx_i)}^2 \notag \\
&\leq \left( \dfrac{2C_2}{\epsilon^2 \omega_d} \right)^2 \int_{H_1(\bzero)} \int_{H_1(\bzero)} J(\abs{\bxi}) J(\abs{\betanew}) \dfrac{1}{\abs{\bxi}} \dfrac{1}{\abs{\betanew}} \notag \\
&\quad  \left( \abs{\be^k(u)(\bx_i + \epsilon \bxi)} + \abs{\be^k(u)(\bx_i )} \right) \left( \abs{\be^k(u)(\bx_i + \epsilon \betanew)} + \abs{\be^k(u)(\bx_i )} \right) d\bxi d\betanew .
\end{align*}
Using inequality $|ab| \leq (\abs{a}^2 + \abs{b}^2)/2$, we get
\begin{align*}
&\left( \abs{\be^k(u)(\bx_i + \epsilon \bxi)} + \abs{\be^k(u)(\bx_i )} \right) \left( \abs{\be^k(u)(\bx_i + \epsilon \betanew)} + \abs{\be^k(u)(\bx_i )} \right) \notag \\
&\leq 3 \left( \abs{\be^k(u)(\bx_i + \epsilon \bxi)}^2 + \abs{\be^k(u)(\bx_i + \epsilon \betanew)}^2 + \abs{\be^k(u)(\bx_i)}^2 \right),
\end{align*}
and
\begin{align*}
&\sum_{\bx_i \in D} h^d \abs{ -\del PD^{\epsilon}(\buhat^k)(\bx_i) + \del PD^{\epsilon}(\butilde^k)(\bx_i)}^2 \notag \\
&\leq \left( \dfrac{2C_2}{\epsilon^2 \omega_d} \right)^2 \int_{H_1(\bzero)} \int_{H_1(\bzero)} J(\abs{\bxi}) J(\abs{\betanew}) \dfrac{1}{\abs{\bxi}} \dfrac{1}{\abs{\betanew}} \notag \\
&\quad \sum_{\bx_i \in D} h^d 3 \left( \abs{\be^k(u)(\bx_i + \epsilon \bxi)}^2 + \abs{\be^k(u)(\bx_i + \epsilon \betanew)}^2 + \abs{\be^k(u)(\bx_i)}^2 \right) d\bxi d\betanew .
\end{align*}
Since $\be^k(u)(\bx) = \sum_{\bx_i \in D} \be^k_i(u) \chi_{U_i}(\bx) $, we have
\begin{align}\label{eq:estimate diff in per force stability}
\sum_{\bx_i \in D} h^d \abs{ -\del PD^{\epsilon}(\buhat^k)(\bx_i) + \del PD^{\epsilon}(\butilde^k)(\bx_i)}^2 &\leq  \dfrac{(6\bar{C})^2}{\epsilon^4} \Ltwonorm{\be^k(u)}{\perdrd}^2.
\end{align}
where $\bar{C}$ is given by \autoref{eq: bar C}. In summary \autoref{eq:estimate diff in per force stability} shows the Lipschitz continuity of the peridynamic force with respect to the $L^2$ norm, see \autoref{eq: lipshitz}, expressed in this context as
\begin{align}\label{eq: lipshitzl2}
\Vert\del PD^{\epsilon}(\buhat^k)(\bx)-  \del PD^{\epsilon}(\butilde^k)\Vert_{L^2(D;\mathbb{R}^d)}\leq\dfrac{(6\bar{C})}{\epsilon^2}\Vert \be^k(u)\Vert_{L^2(D;\mathbb{R}^d)}.
\end{align}

Finally, we substitute above inequality in \autoref{eq:error k ineq 1} to get
\begin{align*}
e^{k+1} &\leq e^k + \Delta t \Ltwonorm{\be^{k+1}(v)}{\perdrd} + \Delta t \tau + \Delta t \dfrac{6\bar{C}}{\epsilon^2}\Ltwonorm{\be^k(u)}{\perdrd}
\end{align*}
We add positive quantity $\Delta t ||e^{k+1}(u)||_{L^2(D;\mathbb{R}^d)} + \Delta t 6\bar{C}/\epsilon^2 ||\be^k(v)||_{L^2(D;\mathbb{R}^d)} $ to the right side of above equation, to get
\begin{align*}
&e^{k+1} \leq ( 1 + \Delta t 6\bar{C}/\epsilon^2)  e^k + \Delta t e^{k+1} + \Delta t \tau \\
\Rightarrow & e^{k+1} \leq \dfrac{( 1 + \Delta t 6\bar{C}/\epsilon^2)}{1 - \Delta t} e^{k} + \dfrac{\Delta t }{1 - \Delta t} \tau.
\end{align*}
We recursively substitute $e^{j}$ on above as follows
\begin{align}
e^{k+1} &\leq \dfrac{( 1 + \Delta t 6\bar{C}/\epsilon^2)}{1 - \Delta t} e^k + \dfrac{\Delta t }{1 - \Delta t} \tau  \notag \\
&\leq \left(\dfrac{( 1 + \Delta t 6\bar{C}/\epsilon^2)}{1 - \Delta t} \right)^2 e^{k-1} + \dfrac{\Delta t }{1 - \Delta t} \tau  \left(1 + \dfrac{( 1 + \Delta t 6\bar{C}/\epsilon^2)}{1 - \Delta t}\right) \notag \\
&\leq ...\notag \\
&\leq \left(\dfrac{( 1 + \Delta t 6\bar{C}/\epsilon^2)}{1 - \Delta t} \right)^{k+1} e^0 + \dfrac{\Delta t }{1 - \Delta t} \tau  \sum_{j=0}^k \left(\dfrac{( 1 + \Delta t 6\bar{C}/\epsilon^2)}{1 - \Delta t} \right)^{k-j}. \label{eq:ek estimnate}
\end{align}

Since $1/(1-\Delta t)= 1 + \Delta t + \Delta t^2 + O(\Delta t^3)$, we have
\begin{align*}
\dfrac{( 1 + \Delta t 6\bar{C}/\epsilon^2)}{1 - \Delta t} &\leq 1 + (1 + 6\bar{C}/\epsilon^2) \Delta t + (1 + 6\bar{C}/\epsilon^2) \Delta t^2 + O(\bar{C}/\epsilon^2) O(\Delta t^3).
\end{align*}
Now, for any $k \leq T /\Delta t$, using identity $(1+ a)^k \leq \exp [ka]$ for $a\leq 0$, we have
\begin{align*}
&\left( \dfrac{ 1 + \Delta t 6\bar{C}/\epsilon^2}{1 - \Delta t} \right)^k \\
&\leq \exp \left[k (1 + 6\bar{C}/\epsilon^2) \Delta t + k(1 + 6\bar{C}/\epsilon^2) \Delta t^2 + k O(\bar{C}/\epsilon^2) O(\Delta t^3) \right] \\
&\leq \exp \left[T (1 + 6\bar{C}/\epsilon^2) + T (1 + 6\bar{C}/\epsilon^2) \Delta t + O(T\bar{C}/\epsilon^2) O(\Delta t^2) \right].
\end{align*}
We write above equation in more compact form as follows
\begin{align*}
&\left( \dfrac{ 1 + \Delta t 6\bar{C}/\epsilon^2}{1 - \Delta t} \right)^k \\
&\leq \exp \left[T (1 + 6\bar{C}/\epsilon^2) (1 + \Delta t + O(\Delta t^2)) \right]. 
\end{align*}

We use above estimate in \autoref{eq:ek estimnate} and get following inequality for $e^k$
\begin{align*}
e^{k+1} &\leq \exp \left[T (1 + 6\bar{C}/\epsilon^2) (1 + \Delta t + O(\Delta t^2)) \right] \left( e^0 + (k+1) \tau \Delta t/(1- \Delta t) \right) \notag \\
&\leq \exp \left[T (1 + 6\bar{C}/\epsilon^2) (1 + \Delta t + O(\Delta t^2)) \right] \left( e^0 + T\tau (1 + \Delta t + O(\Delta t^2) \right).
\end{align*}
where we used the fact that $1/(1-\Delta t) = 1+ \Delta t + O(\Delta t^2)$.

Assuming the error in initial data is zero, i.e. $e^0= 0$, and noting the estimate of $\tau$ in \autoref{eq:estimate tau}, we have
\begin{align*}
&\sup_k e^k \leq \exp \left[T (1 + 6\bar{C}/\epsilon^2) \right] T \tau 
\end{align*}
and we conclude to leading order that
\begin{align}\label{eq: fund est}
\sup_k e^k \leq \exp \left[T (1 + 6\bar{C}/\epsilon^2) \right] T \left[ C_t \Delta t + (C_s/\epsilon^2) h^\gamma \right],
\end{align}
Here the constants $C_t$ and $C_s$ are given by \autoref{eq:const Ct} and \autoref{eq:const Cs}.
 This shows the stability of the numerical scheme. We now address the general one step time discretization.

\subsection{Extension to the implicit schemes}\label{ss:implicit}
Let $\theta \in [0, 1]$ be the parameter which controls the contribution of the implicit and explicit scheme. Let $(\buhat^k, \bvhat^k)$ be the solution of \autoref{eq:finite diff eqn u general} and \autoref{eq:finite diff eqn v general} for given fixed $\theta$.

The forward Euler scheme, backward Euler scheme, and Crank Nicholson scheme correspond to the choices $\theta = 0$, $\theta = 1$, and $\theta = 1/2$ respectively. 

To simplify the equations, we define $\Theta$ acting on discrete set $\{f^k\}_k$ as $\Theta f^k := (1-\theta) f^k + \theta f^{k+1}$. By $\Theta \norm{f^k}$, we mean $(1-\theta) \norm{f^k} + \theta \norm{f^{k+1}}$. Following the same steps as in the case of forward Euler, we write down the equation for $\be^k_i(u) := \buhat^k_i - \butilde^k_i$ and $\be^k_i(v) := \bvhat^k_i - \bvtilde^k_i$ as follows
\begin{align}
\be^{k+1}_i(u) &= \be^k_i(u) + \Delta t \Theta \be^k_i(v) + \Delta t \Theta \tau^k_i(u) \label{eq:finite diff error eqn u general} \\
\be^{k+1}_i(v) &= \be^k_i(v) + \Delta t \Theta \sigma^k_i(u) + \Delta t \Theta \sigma^k_i(v) + \Delta t \Theta \tau^k_i(v) \notag  \\
&\quad  + \Delta t (1-\theta) \left(-\del PD^{\epsilon} (\buhat^k)(\bx_i) + \del PD^{\epsilon}(\butilde^k)(\bx_i) \right) \notag \\
&\quad + \Delta t \theta \left(-\del PD^{\epsilon} (\buhat^{k+1})(\bx_i) + \del PD^{\epsilon}(\butilde^{k+1})(\bx_i) \right). \label{eq:finite diff error eqn v general}
\end{align}
where $\tau^k_i(v), \sigma^k_i(u) ,\sigma^k_i(v)$ are defined in \autoref{eq:consistency error in v}, \autoref{eq:consistency error in u spatial}, and \autoref{eq:consistency error in v spatial} respectively. In this section $\tau^k(u)$ is defined as follows
\begin{align*}
\tau^k_i(u) &:= \dparder{\butilde^k_i}{t} - \dfrac{\butilde^{k+1}_i - \butilde^k_i}{\Delta t}.
\end{align*}

We take the $L^2$ norm of $\be^k(u)(\bx)$ and $\be^k(v)(\bx)$ and for brevity we denote the $L^2$ norm by $||\cdot||$. Recall that $\be^k(u)$ and $\be^k(v)$ are the piecewise constant extension of $\{ \be^k_i(u)\}_i$ and $\{ \be^k_i(v)\}_i$ and we get
\begin{align}
\norm{\be^{k+1}(u)} &\leq \norm{\be^k(u)} + \Delta t \Theta \norm{\be^k(v)} + \Delta t \Theta \norm{\tau^k(u)} \label{eq:error norm be u} \\
\norm{\be^{k+1}(v)} &\leq \norm{\be^k(v)} + \Delta t \left( \Theta \norm{\sigma^k(u)} + \Theta \norm{\sigma^k(v)} + \Theta \norm{\tau^k(v)} \right) \notag \\
&\quad +\Delta t (1-\theta)\left( \sum_{i, \bx_i \in D} h^d \abs{-\del PD^{\epsilon} (\buhat^k)(\bx_i) + \del PD^{\epsilon}(\butilde^k)(\bx_i)}^2 \right)^{1/2} \notag \\
&\quad + \Delta t \theta \left( \sum_{i, \bx_i \in D} h^d \abs{-\del PD^{\epsilon} (\buhat^{k+1})(\bx_i) + \del PD^{\epsilon}(\butilde^{k+1})(\bx_i)}^2 \right)^{1/2}. \label{eq:error norm be v}
\end{align}

From our consistency analysis, we have
\begin{align}
\tau &= \sup_k \left(\Ltwonorm{\tau^k(u)}{\perdrd}  + \Ltwonorm{\tau^k(v)}{\perdrd} \right. \notag \\
&\quad \left. + \Ltwonorm{\sigma^k(u)}{\perdrd} + \Ltwonorm{\sigma^k(v)}{\perdrd}\right) \notag \\
&\leq C_t \Delta t + C_s\dfrac{h^\gamma}{\epsilon^2}. \label{eq:estimate tau 1}
\end{align}
where $C_t$ and $C_s$ are given by \autoref{eq:const Ct} and \autoref{eq:const Cs}. Since $0\leq 1-\theta \leq 1$ and $0\leq \theta \leq 1$ for all $\theta \in [0,1]$, we have
\begin{align*}
\Theta \left(\norm{\tau^k(u)}  + \norm{\tau^k(v)} + \norm{\sigma^k(u)} +\norm{\sigma^k(v)} \right) &\leq 2 \tau.
\end{align*}

\textbf{Crank Nicholson scheme: }If $\theta = 1/2$, and if $\bu,\bv \in \Cthreeintime{[0,T]}{\Cholder{\perdrd}}$, then we can show that
\begin{align*}
\dfrac{1}{2} \tau^k_i(u) + \dfrac{1}{2} \tau^{k+1}_i(u) = \dfrac{(\Delta t)^2}{12} \dfrac{\partial^3 \butilde^{k+1/2}_i}{\partial t^3} + O((\Delta t)^3).
\end{align*}
A similar result holds for $1/2 \tau^k_i(v) + 1/2\tau^{k+1}_i(v)$. Therefore, the consistency error will be bounded by $ \bar{C}_t\Delta t^2 + C_s h^\gamma/\epsilon^2$ with 
\begin{align}\label{eq: timebound}
\bar{C}_t := \frac{1}{12} \sup_{t} || \dfrac{\partial^3 \bu(t)}{\partial t^3} ||_{L^2(D;\bbR^d)} + \frac{1}{12} \sup_{t} || \dfrac{\partial^4 \bu(t)}{\partial t^4} ||_{L^2(D;\bbR^d)} 
\end{align}
and $C_s$ is given by \autoref{eq:const Cs}.

We now estimate \autoref{eq:error norm be v}. Similar to \autoref{eq:estimate diff in per force stability}, we have
\begin{align}
\left( \sum_{i, \bx_i \in D} h^d \abs{-\del PD^{\epsilon} (\buhat^k)(\bx_i) + \del PD^{\epsilon}(\butilde^k)(\bx_i)}^2 \right)^{1/2} &\leq \dfrac{\bar{C}}{\epsilon^2} \norm{\be^k(u)}, \label{eq:estimate diff in per force stability 1} \\
\left( \sum_{i, \bx_i \in D} h^d \abs{-\del PD^{\epsilon} (\buhat^{k+1})(\bx_i) + \del PD^{\epsilon}(\butilde^{k+1})(\bx_i)}^2 \right)^{1/2} &\leq \dfrac{\bar{C}}{\epsilon^2} \norm{\be^{k+1}(u)}, \label{eq:estimate diff in per force stability 2}
\end{align}
where $\bar{C}$ is the constant given by \autoref{eq: bar C}. 
Let $e^k := \norm{\be^k(u)} + \norm{\be^k(v)}$. Adding \autoref{eq:error norm be u} and \autoref{eq:error norm be v} and noting \autoref{eq:estimate diff in per force stability 1}, \autoref{eq:estimate diff in per force stability 2}, and \autoref{eq:estimate tau 1}, we get
\begin{align*}
e^{k+1} &\leq (1+ \Delta t (1- \theta) \dfrac{\bar{C}}{\epsilon^2} ) e^k + \Delta t \theta \dfrac{\bar{C}}{\epsilon^2} e^{k+1} + 2\tau \Delta t, 
\end{align*}
where we assumed $\bar{C}/\epsilon^2 \geq 1$. We further simplify the equation and write
\begin{align}\label{eq:error general time stepping 1}
e^{k+1} &\leq \dfrac{1+ \Delta t (1- \theta) \bar{C}/\epsilon^2}{1 - \Delta t \theta \bar{C}/\epsilon^2}  e^k+ \dfrac{2 }{1 - \Delta t \theta \bar{C}/\epsilon^2} \tau \Delta t,
\end{align}
where we have assumed that $1 - \Delta t \theta \bar{C}/\epsilon^2 > 0 $, i.e. 
\begin{align}\label{eq:assumption delta t}
\Delta t < \dfrac{\epsilon^2}{\bar{C}}=K\epsilon^2.
\end{align}
Thus, for fixed $\epsilon>0$, the error calculation in this section applies when the time step $\Delta t$  satisfies \autoref{eq:assumption delta t}. We now define $a$ and $b$ by
\begin{align*}
a &:= \dfrac{1+ \Delta t (1- \theta) \bar{C}/\epsilon^2}{1 - \Delta t \theta \bar{C}/\epsilon^2} \\
b &:= \dfrac{1}{1 - \Delta t \theta \bar{C}/\epsilon^2}.
\end{align*}
We use the fact that, for $\Delta t$ small, $(1-\alpha \Delta t)^{-1} = 1 + \alpha \Delta t + \alpha^2 (\Delta t)^2 + O((\Delta t)^3)$, to get 
\begin{align*}
b &= 1 + \Delta t \theta \bar{C}/\epsilon^2 + O\left( \left(\Delta t/\epsilon^2 \right)^2 \right) = 1 + O\left( \Delta t/\epsilon^2 \right).
\end{align*}
Now since $\Delta t < \epsilon^2/\bar{C}$, we have
\begin{align}\label{eq:bound on b}
b = O(1).
\end{align}

We have the estimates for  $a$ given by
\begin{align*}a &\leq (1 + \Delta t (1-\theta) \bar{C}/\epsilon^2) (1 + \Delta t \theta \bar{C}/\epsilon^2 + O((\Delta t/\epsilon^2)^2)) \notag \\
&= 1 + \Delta t (\theta + (1-\theta) ) \bar{C}/\epsilon^2 + O((\Delta t/\epsilon^2)^2) \notag \\
&= 1 + \Delta t \bar{C}/\epsilon^2 + O((\Delta t/\epsilon^2)^2).
\end{align*}

Therefore, for any $k \leq T /\Delta t$, we have
\begin{align*}
a^k &\leq \exp \left[ k \Delta t \dfrac{\bar{C}}{\epsilon^2} + k O\left( \left( \dfrac{\Delta t}{\epsilon^2} \right)^2 \right)  \right] \notag \\
&\leq \exp \left[T\bar{C}/\epsilon^2 + O\left(\Delta t \left( \dfrac{1}{\epsilon^2} \right)^2 \right) \right] \notag \\
&\leq \exp [T C/\epsilon^2+O(\frac{1}{\epsilon^2})],
\end{align*}
where we simplified the bound by incorporating \autoref{eq:assumption delta t}. Then, from \autoref{eq:error general time stepping 1}, we get
\begin{align*}
e^{k+1} &\leq a^{k+1} e^0 + 2\tau \left( \Delta t \sum_{j=0}^k a^{j} \right) b.
\end{align*}
From the estimates on $a^k$, we have
\begin{align}\label{eq:bound on sum ak}
\Delta t \sum_{j=0}^k a^{j} &\leq T \exp [T C/\epsilon^2+O(1/\epsilon^2)].
\end{align}
Combining \autoref{eq:bound on b} and \autoref{eq:bound on sum ak}, to get
\begin{align*}
e^{k} &\leq \exp [ T C/\epsilon^2 +O(1/\epsilon^2)] \left( e^0 + 2 T \tau O(1) \right).
\end{align*}
Since $\tau = O(C_t\Delta t + C_sh^\gamma/\epsilon^2)$, we conclude that, for any $\epsilon > 0$ fixed,
\begin{align*}
\sup_{k} e^{k} \leq O(C_t\Delta t + C_sh^\gamma/\epsilon^2).
\end{align*}
Where we assumed $e^0 = 0$. Similarly, for $\theta = 1/2$, we have $\sup_{k} e^{k+1} \leq O(\bar{C}_t(\Delta t)^2 +C_s h^\gamma /\epsilon^2)$. Therefore, the scheme is stable and consistent for any $\theta \in [0,1]$.

\subsection{Stability of the energy for the semi-discrete approximation}
\label{semidiscrete}
We first spatially discretize the peridynamics equation \autoref{eq:per equation}. Let $\{\hat{\bu}_i(t)\}_{i,\bx_i\in D}$ denote the semi-discrete approximate solution which satisfies following, for all $t\in [0,T]$ and $i$ such that $\bx_i\in D$,
\begin{align}\label{eq:fd semi discrete}
\ddot{\hat{\bu}}_i(t) = -\del PD^\epsilon(\hat{\bu}(t))(\bx_i) + \bb_i(t)
\end{align}
where $\hat{\bu}(t)$ is the piecewise constant extension of discrete set $\{\hat{\bu}_i(t) \}_i$ and is defined as
\begin{align}\label{eq:def piecewise ext}
\hat{\bu}(t,\bx) &:= \sum_{i, \bx_i \in D} \hat{\bu}_i(t) \chi_{U_i}(\bx).
\end{align}
The scheme is complemented with the discretized initial conditions $\hat{\bu}_i(0) = \bu_0(\bx_i)$ and $\hat{\bv}_i(0) =\bv_0(\bx_i)$. We apply boundary condition by setting $\hat{\bu}_i(t) = \bzero$ for all $t$ and for all $\bx_i \notin D$.

We have the stability of semi-discrete evolution. 

{\vskip 2mm}
\begin{theorem}
\textbf{Energy stability of the semi-discrete approximation}\\
Let $\{\hat{\bu}_i(t) \}_i$ satisfy \autoref{eq:fd semi discrete} and $\hat{\bu}(t)$ is its piecewise constant extension. Similarly let $\hat{\bb}(t,\bx)$ denote the piecewise constant extension of $\{ \bb(t,\bx_i)\}_{i,\bx_i\in D}$. Then the peridynamic energy $\mathcal{E}^\epsilon$ as defined in \autoref{eq:def energy} satisfies, $\forall t \in [0,T]$,
\begin{align}\label{eq:inequal energy}
\mathcal{E}^\epsilon(\hat{\bu})(t) &\leq \left( \sqrt{\mathcal{E}^\epsilon(\hat{\bu})(0)} + \dfrac{T C}{\epsilon^{3/2}} + \int_0^T ||\hat{\bb}(s)||_{L^2(D;\bbR^d)} ds \right)^2 .
\end{align}
The constant $C$, defined in \autoref{eq:def const stab semi fd}, is independent of $\epsilon$ and $h$. 
\end{theorem}
{\vskip 2mm}

\begin{proof}
We multiply \autoref{eq:fd semi discrete} by $\chi_{U_i}(\bx)$ and sum over $i$ and use definition of piecewise constant extension in \autoref{eq:def piecewise ext} to get
\begin{align*}
\ddot{\hat{\bu}}(t,\bx) &= -\del \hat{PD^\epsilon}(\hat{\bu}(t))(\bx) + \hat{\bb}(t,\bx) \notag \\
&= -\del PD^\epsilon(\hat{\bu}(t))(\bx) + \hat{\bb}(t,\bx) \notag \\
&\quad + (-\del \hat{PD^\epsilon}(\hat{\bu}(t))(\bx) + \del PD^\epsilon(\hat{\bu}(t))(\bx))
\end{align*}
where $-\del \hat{PD^\epsilon}(\hat{\bu}(t))(\bx)$ and $\hat{\bb}(t,\bx)$ are given by
\begin{align*}
-\del \hat{PD^\epsilon}(\hat{\bu}(t))(\bx) &=  \sum_{i,\bx_i \in D} (-\del PD^\epsilon(\hat{\bu}(t))(\bx_i)) \chi_{U_i}(\bx) \notag \\
\hat{\bb}(t,\bx) &= \sum_{i, \bx_i\in D} \bb(t,\bx_i) \chi_{U_i}(\bx).
\end{align*}

We define set as follows
\begin{align}\label{eq:def sigma fd semi}
\sigma(t,\bx) := -\del \hat{PD^\epsilon}(\hat{\bu}(t))(\bx) + \del PD^\epsilon(\hat{\bu}(t))(\bx).
\end{align}
We use following result which we will show after few steps
\begin{align}\label{eq:bd on sigma fd semi}
||\sigma(t)||_{L^2(D;\bbR^d)} \leq \dfrac{C}{\epsilon^{3/2}}.
\end{align}

We then have
\begin{align}\label{eq:fd semi discrete 2}
\ddot{\hat{\bu}}(t,\bx) &= -\del PD^\epsilon(\hat{\bu}(t))(\bx) + \hat{\bb}(t,\bx) + \sigma(t,\bx).
\end{align}
Multiplying above with $\dot{\hat{\bu}}(t)$ and integrating over $D$ to get
\begin{align*}
(\ddot{\hat{\bu}}(t),\dot{\hat{\bu}}(t))  &= (-\del PD^\epsilon(\hat{\bu}(t)), \dot{\hat{\bu}}(t)) \notag \\
&\quad + (\hat{\bb}(t), \dot{\hat{\bu}}(t)) + (\sigma(t),\dot{\hat{\bu}}(t)). 
\end{align*}
Consider energy $\mathcal{E}^\epsilon(\hat{\bu})(t)$ given by \autoref{eq:def energy} and note the identity \autoref{eq:energy relat}, to have
\begin{align*}
\dfrac{d}{dt} \mathcal{E}^\epsilon(\hat{\bu})(t) &= (\hat{\bb}(t), \dot{\hat{\bu}}(t)) + (\sigma(t),\dot{\hat{\bu}}(t)) \notag \\
&\leq \left( ||\hat{\bb}(t)||_{L^2(D;\bbR^d)}  + ||\sigma(t)||_{L^2(D;\bbR^d)} \right) ||\dot{\hat{\bu}}(t)||_{L^2(D;\bbR^d)},
\end{align*}
where we used H\"older inequality in last step. Since $PD^\epsilon(\bu)$ is positive for any $\bu$, we have
\begin{align*}
||\dot{\hat{\bu}}(t)|| &\leq 2 \sqrt{\dfrac{1}{2}||\dot{\hat{\bu}}(t)||_{L^2(D;\bbR^d)}^2 + PD^\epsilon(\hat{\bu}(t)) } = 2 \sqrt{\mathcal{E}^\epsilon(\hat{\bu})(t)}.
\end{align*}
Using above, we get
\begin{align*}
\dfrac{1}{2}\dfrac{d}{dt} \mathcal{E}^\epsilon(\hat{\bu})(t) &\leq \left( ||\hat{\bb}(t)||_{L^2(D;\bbR^d)}  + ||\sigma(t)||_{L^2(D;\bbR^d)} \right) \sqrt{\mathcal{E}^\epsilon(\hat{\bu})(t)}.
\end{align*}
Let $\delta > 0$ is some arbitrary but fixed real number and let $A(t) = \delta + \mathcal{E}^\epsilon(\hat{\bu})(t)$. Then
\begin{align*}
\dfrac{1}{2}\dfrac{d}{dt} A(t) &\leq \left( ||\hat{\bb}(t)||_{L^2(D;\bbR^d)}  + ||\sigma(t)||_{L^2(D;\bbR^d)} \right) \sqrt{A(t)}.
\end{align*}
Using the fact that $\frac{1}{\sqrt{A(t)}} \frac{d}{dt}A(t) = 2 \frac{d}{dt} \sqrt{A(t)}$, we have
\begin{align*}
\sqrt{A(t)} &\leq \sqrt{A(0)} + \int_0^t \left( ||\hat{\bb}(s)||_{L^2(D;\bbR^d)}  + ||\sigma(s)||_{L^2(D;\bbR^d)} \right) ds \\
&\leq \sqrt{A(0)} + \dfrac{T C}{\epsilon^{3/2}} + \int_0^T ||\hat{\bb}(s)||_{L^2(D;\bbR^d)} ds.
\end{align*}
where we used bound on $||\sigma(s)||_{L^2(D;\bbR^d)}$ from \autoref{eq:bd on sigma fd semi}. Noting that $\delta > 0$ is arbitrary, we send it to zero to get
\begin{align*}
\sqrt{\mathcal{E}^\epsilon(\hat{\bu})(t)} &\leq \sqrt{\mathcal{E}^\epsilon(\hat{\bu})(0)} + \dfrac{T C}{\epsilon^{3/2}} + \int_0^T ||\hat{\bb}(s)|| ds,
\end{align*}
and \autoref{eq:inequal energy} follows by taking square of above equation.

It remains to show \autoref{eq:bd on sigma fd semi}. To simplify the calculations, we use following notations: let $\bxi \in H_1(\bzero)$ and let
\begin{align*}
& s_{\bxi} = \epsilon |\bxi|, e_{\bxi} = \dfrac{\bxi}{|\bxi|}, \bar{\omega}(\bx) = \omega(\bx) \omega(\bx+\epsilon\bxi), \notag \\
& S_{\bxi}(\bx) = \dfrac{\hat{\bu}(t,\bx+\epsilon \bxi) - \hat{\bu}(t,\bx)}{s_{\bxi}} \cdot e_{\bxi}.
\end{align*}
With above notations and using expression of $-\del PD^\epsilon$ from \autoref{eq:peri force use}, we have for $\bx\in U_i$
\begin{align}
&|\sigma(t,\bx)| = \left\vert -\del PD^\epsilon(\hat{\bu}(t))(\bx_i) + \del PD^\epsilon(\hat{\bu}(t))(\bx) \right\vert \notag \\
&= \left\vert \dfrac{2}{\epsilon \omega_d} \int_{H_1(\bzero)} \dfrac{J(|\bxi|)}{\sqrt{s_{\bxi}}} \left(\bar{\omega}(\bx_i) F'_1(\sqrt{s_{\bxi}} S_{\bxi}(\bx_i)) - \bar{\omega}(\bx) F'_1(\sqrt{s_{\bxi}} S_{\bxi}(\bx)) \right) e_{\bxi} d\bxi \right\vert \notag \\
&\leq \dfrac{2}{\epsilon \omega_d} \int_{H_1(\bzero)} \dfrac{J(|\bxi|)}{\sqrt{s_{\bxi}}} \left\vert \bar{\omega}(\bx_i) F'_1(\sqrt{s_{\bxi}} S_{\bxi}(\bx_i)) - \bar{\omega}(\bx) F'_1(\sqrt{s_{\bxi}} S_{\bxi}(\bx)) \right\vert d\bxi \notag \\
&\leq \dfrac{2}{\epsilon \omega_d} \int_{H_1(\bzero)} \dfrac{J(|\bxi|)}{\sqrt{s_{\bxi}}} \left( \left\vert \bar{\omega}(\bx_i) F'_1(\sqrt{s_{\bxi}} S_{\bxi}(\bx_i)) \right\vert + \left\vert  \bar{\omega}(\bx) F'_1(\sqrt{s_{\bxi}} S_{\bxi}(\bx)) \right\vert \right) d\bxi.
\end{align}
Using the fact that $0\leq \omega(\bx) \leq 1$ and $|F_1'(r)|\leq C_1$, where $C_1$ is $\sup_r |F_1'(r)|$, we get
\begin{align*}
|\sigma(t,\bx)| &\leq \dfrac{4 C_1 \bar{J}_{1/2}}{\epsilon^{3/2}}.
\end{align*}
where $\bar{J}_{1/2} = (1/\omega_d)\int_{H_1(\bzero)} J(|\bxi|) |\bxi|^{-1/2} d\bxi$.

Taking the $L^2$ norm of $\sigma(t,\bx)$, we get
\begin{align*}
||\sigma(t)||_{L^2(D;\bbR^d)}^2 &= \sum_{i,\bx_i\in D} \int_{U_i} |\sigma(t,\bx)|^2 d\bx \leq \left(\dfrac{4 C_1 \bar{J}_{1/2}}{\epsilon^{3/2}}\right)^2 \sum_{i,\bx_i\in D} \int_{U_i} d\bx
\end{align*}
thus
\begin{align*}
||\sigma(t)||_{L^2(D;\bbR^d)} \leq \dfrac{4 C_1 \bar{J}_{1/2}\sqrt{|D|}}{\epsilon^{3/2}} = \dfrac{C}{\epsilon^{3/2}}
\end{align*}
where
\begin{align}\label{eq:def const stab semi fd}
C := 4 C_1 \bar{J}_{1/2}\sqrt{|D|}.
\end{align}

This completes the proof.
\end{proof}

\subsection{Local instability under radial perturbations}
\label{ss:proof localstab}
We observe that both explicit and implicit schemes treated in previous sections show that any increase in local truncation error is controlled at each time step. From the proofs above (and the general approximation theory for ODE), this control is adequate to establish convergence rates as $\Delta t \rightarrow 0$. We now comment on a source of error that can grow with time steps in regions where the strain is large and peridynamic bonds between material points begin to soften.

We examine the Jacobian matrix of the peridynamic system associated by perturbing about a displacement field and seek to understand the stability of the perturbation. Suppose the solution is near the displacement field $\overline{\bu}(\bx)$ and let $\bs(t,\bx)=\bu(t,\bx)-\overline{\bu}(\bx)$ be the perturbation. 
We write the associated strain as $S(\by,\bx;\overline{\bu})$ and $S(\by,\bx;\bs)$.
Expanding the peridynamic force in Taylor series about $\overline{\bu}$ assuming $\bs$ is small gives 
\begin{align*}
\partial_{tt} \bs(t,\bx)&=\frac{2}{V_d}\left\{\int_{\mathcal{H}_\epsilon(\bx)}\partial^2_{{S}}\mathcal{W}^\epsilon({S}(\by,\bx;\overline{\bu}))S(\by,\bx;\bs)\frac{\by-\bx}{|\by-\bx|}d\by\right\} \notag \\
&\quad -\del PD^{\epsilon}(\overline{\bu})(\bx) +\bb(t,\bx) + O(|\bs|^2),
\end{align*}
where $\mathcal{W}^\epsilon(S, \by - \bx) = W^\epsilon(S,\by - \bx)/|\by - \bx|$ and $W^\epsilon$ is given by \autoref{eq:per pot}.

To recover a local stability formula in terms of a spectral radius we consider local radial perturbations $\bs$ with spatially constant strain $S(\by,\bx;\bs)$ of the form $S(\by,\bx;\bs)=-\delta(t){\bolds{\mu}}\cdot\be$ where $\bolds{\mu}$ is in $\mathbb{R}^d$ and $\bs$ has radial variation about $\bx$ with $\bs(\by)=\delta(t){\bolds{\mu}}(1-|\by-\bx|)$. This delivers the local ODE
\begin{align*}
\delta''(t)\bolds{\mu}=A\delta(t){\bolds{\mu}}+b
\end{align*}
where the stability matrix $A$ is selfadjoint and given by
\begin{align}\label{eq:A}
A=-\frac{2}{V_d}\left\{\int_{\mathcal{H}_\epsilon(\bx)}\partial^2_{{S}}\mathcal{W}^\epsilon({S}(\by,\bx;\overline{\bu}))\frac{\by-\bx}{|\by-\bx|}\otimes\frac{\by-\bx}{|\by-\bx|} d\by\right\},
\end{align}
and
\begin{align*}
b=-\del PD^{\epsilon}(\overline{\bu})(\bx) +\bb(t,\bx) + O(|\bs|^2).
\end{align*}
A stability criterion for the perturbation is obtained on analyzing the linear system $\delta''(t)\bolds{\mu}=A\delta(t){\bolds{\mu}}$.  Writing it as a 1st order system gives
\begin{align*}
\delta_1'(t)\bolds{\mu}&=\delta_2(t)\bolds{\mu} \notag \\
\delta_2'(t)\bolds{\mu}&=A\delta_1(t)\bolds{\mu}
\end{align*}
where $\bolds{\mu}$ is a vector in $\mathbb{R}^d$.
The eigenvalues of $A$ are real and denoted by $\lambda_i$, $i=1,\ldots,d$ and the associated eigenvectors are denoted by $\bv^i$.
Choosing $\bolds{\mu}=\bv^i$ gives
\begin{align*}
\delta_1'(t)\bv^i&=\delta_2(t)\bv^i\notag \\
\delta_2'(t)\bv^i&={\lambda_i}\delta_1(t){\bv^i}.
\end{align*}
Applying the Forward Euler method to this system gives the discrete iterative system 
\begin{align*}
\delta_1^{k+1} &= \delta_1^k+\Delta t \delta_2^k \notag \\
\delta_2^{k+1} &= {\lambda_i}\Delta t \delta_1^k +\delta_2^k.
\end{align*}
The spectral radius of the matrix associated with this iteration is 
\begin{align*}
\rho=\max_{i=1,\ldots,d}|1\pm \Delta t \sqrt{\lambda_i}|.
\end{align*}
It is easy to see that the spectral radius is larger than $1$ for any choice of $\lambda_i$ and we conclude local instability for the forward Euler scheme under radial perturbation.

For the implicit scheme given by backward Euler we get the discrete iterative system
\begin{align*}
\delta_1^{k} &= \delta_1^{k+1}-\Delta t \delta_2^{k+1} \notag \\
\delta_2^{k} &= {-\lambda_i}\Delta t \delta_1^{k+1} +\delta_2^{k+1}.
\end{align*}
and
\begin{equation*}
  \left[\begin{array}{c}
   \delta_1^{k+1} \\    
   \delta_2^{k+1}
    \end{array}\right ]=
  \left[\begin{array}{cc}
    1 &-\Delta t\\        
    -\Delta t  \lambda_i&1
  \end{array}\right]^{-1}
  \left[\begin{array}{c}
   \delta_1^{k} \\    
   \delta_2^{k}
  \end{array}\right].
  \label{eq:baclwardmatrix}
\end{equation*}
The spectral radius for the iteration matrix is
\begin{align*}\label{eq:diccretestabbackwards}
\rho=\max_{i=1,\ldots,d}|\frac{1}{\theta}\pm \frac{\Delta t \sqrt{\lambda_i}|}{|\theta|}|,
\end{align*}
where $\theta=1-\lambda_i (\Delta t)^2$. If we suppose that the stability matrix $A$ is not negative definite and there is a  $\lambda_j>0$, then the spectral radius is larger than one, i.e, 
\begin{align}\label{eq:discretestabbackwardsexplicit}
1<\frac{|1+ \sqrt{\lambda_j}\Delta t|}{|1-\lambda_j(\Delta t)^2|}\leq\rho.
\end{align}
Thus it follows from \autoref{eq:discretestabbackwardsexplicit}
that we can have local instability of the backward Euler scheme for radial perturbations. 
Inspection of \autoref{eq:A} shows the sign of the eigenvalues of the matrix $A$ depend explicitly on the sign of $\partial^2_{{S}}\mathcal{W}^\epsilon({S}(\by,\bx;\overline{\bu})$. It is shown in  \cite{CMPer-Lipton} that 
\begin{align}
\partial^2_{{S}}\mathcal{W}^\epsilon({S}(\by,\bx;\overline{\bu})>0 \hbox{  for  } |{S}(\by,\bx;\overline{\bu})|<S_c \label{eq:stabilityb}\\
\partial^2_{{S}}\mathcal{W}^\epsilon({S}(\by,\bx;\overline{\bu})<0 \hbox{  for  } |{S}(\by,\bx;\overline{\bu})|>S_c \label{eq:nostability}.
\end{align}
From the model we see that bonds are loosing stiffness when $|{S}(\by,\bx;\overline{\bu})|>S_c $ and the points for which $ A$ is non negative definite correspond to points where \autoref{eq:nostability} holds for a preponderance of bonds inside the horizon. 
We conclude noting that both explicit and implicit schemes treated in previous sections have demonstrated convergence rates  $O((C_t\Delta t + C_sh^\gamma/\epsilon^2))$ as $\Delta t \rightarrow 0$. However, the results of this section show that the error can grow with time for this type of radial perturbation.

\section{Lipschitz continuity in H\"older norm and existence of a solution}
\label{s:proof existence}
In this section, we prove \autoref{prop:lipschitz}, \autoref{thm:local existence}, and \autoref{thm:existence over finite time domain}.

\subsection{Proof of Proposition 1}
Let $I = [0,T]$ be the time domain and $X=\Cholderz{\perdrd}\times\Cholderz{\perdrd}$. Recall that $F^\epsilon(y,t) = (F^\epsilon_1(y,t), F^\epsilon_2(y,t))$, where $F^\epsilon_1(y,t) = y^2$ and $F^\epsilon_2(y,t) = -\del PD^\epsilon(y^1) + \bb(t)$. Given $t\in I$ and $y=(y^1,y^2), z=(z^1, z^2)\in X$, we have
\begin{align}\label{eq:norm of F in X}
&\normX{F^\epsilon(y,t) - F^\epsilon(z,t)}{X} \notag \\
&\leq \Choldernorm{y^2 - z^2}{\perdrd} + \Choldernorm{-\del PD^{\epsilon}(y^1) + \del PD^\epsilon(z^1)}{\perdrd}.
\end{align}

Therefore, to prove the \autoref{eq:lipschitz property of F}, we only need to analyze the second term in above inequality. Let $\bu,\bv \in \Cholderz{\perdrd}$, then we have
\begin{align}\label{eq:lipschitz norm of per force}
&\Choldernorm{-\del PD^{\epsilon}(\bu) - (-\del PD^{\epsilon}(\bv))}{\perdrd} \notag \\
& = \sup_{\bx\in D} \abs{-\del PD^{\epsilon}(\bu)(\bx) - (-\del PD^{\epsilon}(\bv)(\bx))}  \notag \\
&\quad + \sup_{\substack{\bx\neq \by,\\
\bx,\by \in D}} \dfrac{\abs{(-\del PD^{\epsilon}(\bu) + \del PD^{\epsilon}(\bv))(\bx) - (-\del PD^{\epsilon}(\bu) + \del PD^{\epsilon}(\bv))(\by)}}{\abs{\bx - \by}^{\gamma}}.
\end{align}
Note that the force $-\del PD^{\epsilon}(\bu)(\bx)$ can be written as follows
\begin{align*}
&- \del PD^{\epsilon}(\bu)(\bx) \notag \\
&= \dfrac{4}{\epsilon^{d+1} \omega_d} \int_{H_\epsilon(\bx)} \omega(\bx)\omega(\by)J(\dfrac{\abs{\by - \bx}}{\epsilon}) f'(\abs{\by - \bx} S(\by, \bx; \bu)^2) S(\by,\bx;\bu) \dfrac{\by - \bx}{\abs{\by - \bx}} d\by \notag \\
&= \dfrac{4}{\epsilon \omega_d} \int_{H_1(\bzero)} \omega(\bx)\omega(\bx+\epsilon \bxi)J(\abs{\bxi}) f'(\epsilon \abs{\bxi} S(\bx+\epsilon \bxi, \bx; \bu)^2) S(\bx+ \epsilon \bxi,\bx;\bu) \dfrac{\bxi}{\abs{\bxi}} d\bxi.
\end{align*}
where we substituted $\partial_S W^{\epsilon}$ using \autoref{eq:per pot}. In second step, we introduced the change in variable $\by = \bx + \epsilon \bxi$. 

Let $F_1: \bbR \to \bbR$ be defined as $F_1(S) = f(S^2)$. Then $F'_1(S) = f'(S^2) 2S$. Using the definition of $F_1$, we have
\begin{align*}
2 S f'(\epsilon \abs{\bxi} S^2) = \dfrac{F'_1(\sqrt{\epsilon \abs{\bxi}} S)}{\sqrt{\epsilon \abs{\bxi}}}.
\end{align*} 

\begin{figure}
\centering
\includegraphics[scale=0.3]{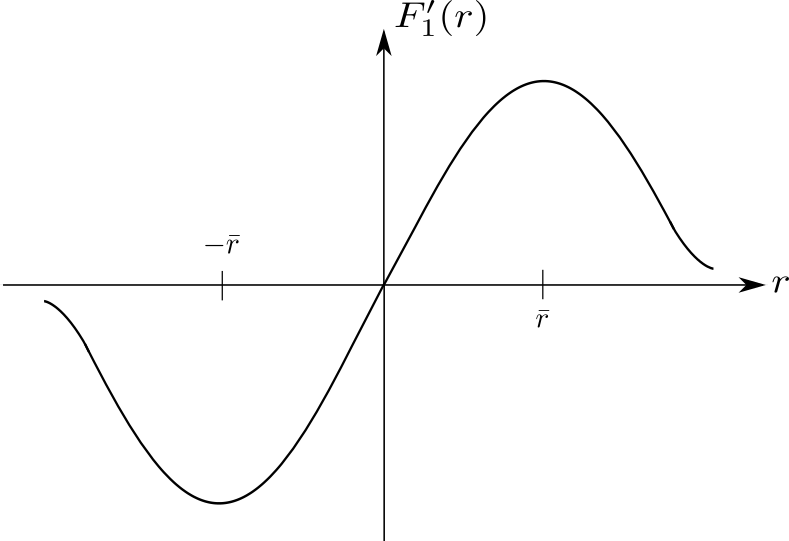}
\caption{Generic plot of $F'_1(r)$. $|F'_1(r)|$ is bounded by $\abs{F'_1(\bar{r})}$.}
 \label{fig:first der per pot one}
\end{figure}

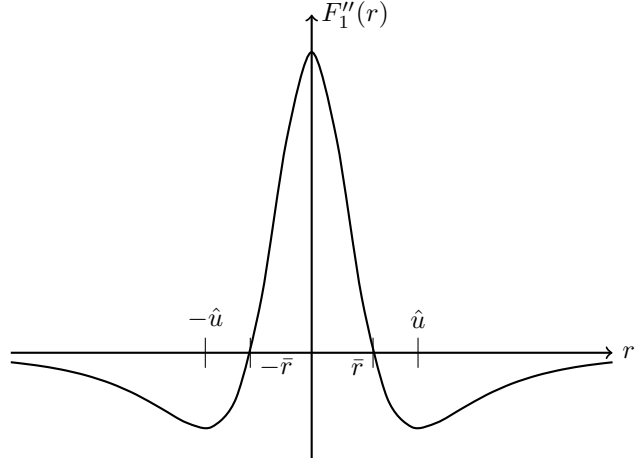
\begin{figure}
\centering
\begin{tikzpicture}
      \draw[->,thick] (-4,0) -- (4,0) node[right] {$r$};
      \draw[->,thick] (0,-1.4) -- (0,4.5) node[right] {$F''_1(r)$};
      \draw[scale=1.0,domain=-4:4,smooth,variable=\x,thick] plot ({\x},{\dderperpotFx(\x)});
      \draw [-] (0.816,0.2) -- (0.816,-0.2);
      \draw [-] (-0.816,0.2) -- (-0.816,-0.2);
      \draw [-] (1.414,0.2) -- (1.414,-0.2);
      \draw [-] (-1.414,0.2) -- (-1.414,-0.2);
      \node [left] at (0.816,-0.2) {$\bar{r}$};
      \node [right] at (-0.816,-0.2) {$-\bar{r}$};
      \node [above] at (1.414,0.2) {$\hat{u}$};
      \node [above] at (-1.414,0.2) {$-\hat{u}$};         
\end{tikzpicture}

\caption{Generic plot of $F''_1(r)$. At $\pm\bar{r}$, $F''_1(r) = 0$. At $\pm\hat{u}$, $F'''_1(r) = 0$.}
 \label{fig:second der per pot}
\end{figure}

\begin{figure}
\centering
\begin{tikzpicture}
      \draw[->,thick] (-5,0) -- (5,0) node[right] {$r$};
      \draw[->,thick] (0,-1) -- (0,1.3) node[above] {$F'''_1(r)$};
      \draw[scale=1.0,domain=-5:5,smooth,variable=\x,thick] plot ({\x},{\ddderperpotFx(\x)});
      \draw[-] (0.8, 0.1) -- (0.8, -0.1);
      \draw[-] (-0.8, 0.1) -- (-0.8, -0.1);
       \node [above] at (0.8, 0.1) {$\bar{u}_2$};
       \node [below] at (-0.8,-0.1) {$-\bar{u}_2$};
      \draw[-] (1.64, 0.1) -- (1.64, -0.1);
      \draw[-] (-1.64, 0.1) -- (-1.64, -0.1);
       \node [right] at (1.64, -0.2) {$\hat{u}$};
       \node [left] at (-1.64, 0.2) {$-\hat{u}$};
      \draw[-] (2.9, 0.1) -- (2.9, -0.1);
      \draw[-] (-2.9, 0.1) -- (-2.9, -0.1);
       \node [below] at (2.9, -0.1) {$\tilde{u}_2$};
       \node [above] at (-2.9, 0.1) {$-\tilde{u}_2$};
\end{tikzpicture}
\caption{Generic plot of $F'''_1(r)$. At $\pm \bar{u}_2$ and $\pm \tilde{u}_2$, $F''''_1 = 0$.}
\label{fig:third der per pot}
\end{figure}
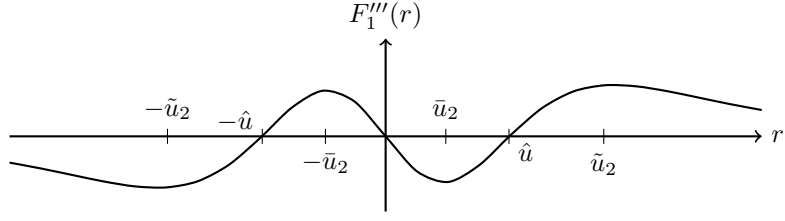

Because $f$ is assumed to be positive, smooth, and concave, and is bounded far away, we have following bound on derivatives of $F_1$
\begin{align}
\sup_{r} \abs{F'_1(r)} &= F'_1(\bar{r}) =: C_1 \label{C1}\\
\sup_{r} \abs{F''_1(r)} &= \max \{ F''_1(0), F''_1(\hat{u}) \} =: C_2 \label{C2}\\
\sup_{r} \abs{F'''_1(r)} &= \max \{F'''_1(\bar{u}_2), F'''_1(\tilde{u}_2) \} =:C_3.\label{C3}
\end{align}
where $\bar{r}$ is the inflection point of $f(r^2)$, i.e. $F''_1(\bar{r}) = 0$. $\{ 0, \hat{u} \}$ are the maxima of $F''_1(r)$. $\{\bar{u}, \tilde{u} \}$ are the maxima of $F'''_1(r)$. By chain rule and by considering the assumption on $f$, we can show that $\bar{r}, \hat{u}, \bar{u}_2, \tilde{u}_2$ exists and the $C_1, C_2, C_3$ are bounded. \autoref{fig:first der per pot one}, \autoref{fig:second der per pot}, and \autoref{fig:third der per pot} shows the generic graphs of $F'_1(r)$, $F''_1(r)$, and $F'''_1(r)$ respectively.

The nonlocal force $-\del PD^\epsilon$ can be written as
\begin{align}\label{eq:peri force use}
&- \del PD^{\epsilon}(\bu)(\bx) \notag\\
&= \dfrac{2}{\epsilon \omega_d} \int_{H_1(\bzero)}  \omega(\bx)\omega(\bx+\epsilon \bxi)J(\abs{\bxi}) F'_1(\sqrt{\epsilon \abs{\bxi}} S(\bx+\epsilon \bxi, \bx; \bu)) \dfrac{1}{\sqrt{\epsilon \abs{\bxi}}} \dfrac{\bxi}{\abs{\bxi}} d\bxi.
\end{align}

To simplify the calculations, we use following notation
\begin{align*}
\bar{\bu}(\bx) &:= \bu(\bx+\epsilon \bxi) - \bu(\bx), \\
\bar{\bu}(\by) &:= \bu(\by+\epsilon \bxi) - \bu(\by), \\
(\bu - \bv)(\bx) &:= \bu(\bx) - \bv(\bx),
\end{align*}
and $\overline{(\bu - \bv)}(\bx)$ is defined similar to $\bar{\bu}(\bx)$. Also, let
\begin{align*}
\quad s = \epsilon \abs{\bxi}, \quad \be = \dfrac{\bxi}{\abs{\bxi}}.
\end{align*}

In what follows, we will come across the integral of type $\int_{H_1(\bzero)} J(\abs{\bxi}) \abs{\bxi}^{-\alpha} d\bxi$. Recall that $0\leq J(\abs{\bxi}) \leq M$ for all $\bxi\in H_1(\bzero)$ and $J(\abs{\bxi}) = 0$ for $\bxi \notin H_1(\bzero)$. Therefore, let
\begin{align}
\label{Jbar}
\bar{J}_\alpha &:= \dfrac{1}{\omega_d} \int_{H_1(\bzero)} J(\abs{\bxi}) \abs{\bxi}^{-\alpha} d\bxi.
\end{align}
With notations above, we note that $S(\bx+ \epsilon \bxi, \bx; \bu) = \bar{\bu}(\bx)\cdot \be/s$. $-\del PD^\epsilon$ can be written as
\begin{align}\label{eq:del pd simplified}
- \del PD^{\epsilon}(\bu)(\bx) &= \dfrac{2}{\epsilon \omega_d} \int_{H_1(\bzero)}  \omega(\bx)\omega(\bx+\epsilon \bxi)J(\abs{\bxi}) F'_1(\bar{\bu}(\bx) \cdot \be/\sqrt{s} ) \dfrac{1}{\sqrt{s}} \be d\bxi.
\end{align}

We first estimate the term $\abs{-\del PD^{\epsilon}(\bu)(\bx) - (-\del PD^{\epsilon}(\bv)(\bx))}$ in \autoref{eq:lipschitz norm of per force}.  
\begin{align}
&\abs{-\del PD^{\epsilon}(\bu)(\bx) - (-\del PD^{\epsilon}(\bv)(\bx))}\notag\\
&\leq\abs{ \dfrac{2}{\epsilon \omega_d} \int_{H_1(\bzero)} \omega(\bx)\omega(\bx+\epsilon \bxi) J(\abs{\bxi}) \dfrac{\left( F'_1(\bar{\bu}(\bx) \cdot \be/\sqrt{s} ) - F'_1(\bar{\bv}(\bx) \cdot \be/\sqrt{s} ) \right) }{\sqrt{s}}\be d\bxi } \notag \\
&\leq \abs{ \dfrac{2}{\epsilon \omega_d} \int_{H_1(\bzero)} J(\abs{\bxi}) \dfrac{1}{\sqrt{s}}  \abs{ F'_1(\bar{\bu}(\bx) \cdot \be/\sqrt{s} ) - F'_1(\bar{\bv}(\bx) \cdot \be/\sqrt{s} ) } d\bxi } \notag \\ 
&\leq \sup_{r} \abs{F''_1(r)} \abs{ \dfrac{2}{\epsilon \omega_d} \int_{H_1(\bzero)} J(\abs{\bxi}) \dfrac{1}{\sqrt{s}}  \abs{ \bar{\bu}(\bx) \cdot \be/\sqrt{s} - \bar{\bv}(\bx) \cdot \be/\sqrt{s} } d\bxi } \notag \\ 
&\leq \dfrac{2C_2 }{\epsilon \omega_d} \abs{\int_{H_1(\bzero)} J(\abs{\bxi}) \dfrac{\abs{\bubar({\bx}) - \bvbar({\bx})}}{\epsilon \abs{\bxi}} d\bxi }. \label{eq:diff force bound 1}
\end{align}
Here we have used the fact that $|\omega(\bx)|\leq1$ and for a vector $\be$ such that $\abs{\be} = 1$, $\abs{\ba \cdot \be} \leq \abs{\ba}$ holds and $\abs{\alpha \be} \leq \abs{\alpha}$ holds for all $\ba \in \bbR^d, \alpha \in \bbR$. Using the fact that $\bu,\bv \in \Cholderz{\perdrd}$, we have
\begin{align*}
\dfrac{\abs{\bar{\bu}(\bx) - \bar{\bv}(\bx)}}{s} &= \dfrac{\abs{(\bu - \bv)(\bx+ \epsilon \bxi)- (\bu - \bv)(\bx)}}{(\epsilon \abs{\bxi})^{\gamma}} \dfrac{1}{(\epsilon \abs{\bxi})^{1-\gamma}} \\
&\leq \Choldernorm{\bu-\bv}{\perdrd} \dfrac{1}{(\epsilon \abs{\bxi})^{1-\gamma}}.
\end{align*}
Substituting the estimate given above, we get
\begin{align}\label{eq:estimate lipschitz part 1}
\abs{-\del PD^{\epsilon}(\bu)(\bx) - (-\del PD^{\epsilon}(\bv)(\bx))} &\leq \dfrac{2C_2 \bar{J}_{1-\gamma}}{\epsilon^{2-\gamma}} \Choldernorm{\bu-\bv}{\perdrd},
\end{align}
where $C_2$ is given by \autoref{C2} and $\bar{J}_{1-\gamma}$ is given by \autoref{Jbar}.

We now estimate the second term in \autoref{eq:lipschitz norm of per force}. To simplify notation we write $\tilde{\omega}(\bx,\bxi)=\omega(\bx)\omega(\bx+\epsilon\bxi)$ and with the help of \autoref{eq:del pd simplified}, we get
\begin{align}
&\dfrac{1}{\abs{\bx - \by}^{\gamma}} \abs{(-\del PD^{\epsilon}(\bu) + \del PD^{\epsilon}(\bv))(\bx) - (-\del PD^{\epsilon}(\bu) + \del PD^{\epsilon}(\bv))(\by)} \notag \\
&= \dfrac{1}{\abs{\bx - \by}^{\gamma}} | \dfrac{2}{\epsilon\omega_d} \int_{H_1(0)} J(\abs{\bxi})  \dfrac{1}{\sqrt{s}}\times \left( \tilde{\omega}(\bx,\bxi)(F'_1(\dfrac{\bubar(\bx)\cdot \be}{\sqrt{s}}) - F'_1(\dfrac{\bvbar(\bx))\cdot \be}{\sqrt{s}}))\right.\notag \\
& \: \left. - \tilde{\omega}(\by,\bxi)(F'_1(\dfrac{\bubar(\by)\cdot \be}{\sqrt{s}}) - F'_1(\dfrac{\bvbar(\by))\cdot \be}{\sqrt{s}}) \right) \be d\bxi | \notag \\
&\leq \dfrac{1}{\abs{\bx - \by}^{\gamma}} | \dfrac{2}{\epsilon\omega_d} \int_{H_1(0)} J(\abs{\bxi}) \dfrac{1}{\sqrt{s}}\times \notag \\
&\quad | \tilde{\omega}(\bx,\bxi)(F'_1(\dfrac{\bubar(\bx)\cdot \be}{\sqrt{s}}) - F'_1(\dfrac{\bvbar(\bx)\cdot \be}{\sqrt{s}})) - \tilde{\omega}(\by,\bxi)(F'_1(\dfrac{\bubar(\by)\cdot e}{\sqrt{s}}) - F'_1(\dfrac{\bvbar(\by)\cdot e}{\sqrt{s}})) |  d\bxi.\notag\\ \label{eq:H def}
\end{align}

We analyze the integrand in above equation. We let $H$ be defined by
\begin{align*}
H &:= \dfrac{| \tilde{\omega}(\bx,\bxi)(F'_1(\dfrac{\bubar(\bx)\cdot \be}{\sqrt{s}}) - F'_1(\dfrac{\bvbar(\bx)\cdot \be}{\sqrt{s}})) - \tilde{\omega}(\by,\bxi)(F'_1(\dfrac{\bubar(\by)\cdot e}{\sqrt{s}}) - F'_1(\dfrac{\bvbar(\by)\cdot e}{\sqrt{s}})) |}{\abs{\bx - \by}^\gamma}.
\end{align*}
Let $\br : [0,1] \times D \to \bbR^d$ be defined as 
\begin{align*}
\br(l,\bx) = \bvbar(\bx) + l (\bubar(\bx) - \bvbar(\bx)).
\end{align*}

Note $\partial \br(l,\bx)/ \partial l = \bubar(\bx) - \bvbar(\bx)$. Using $\br(l,\bx)$, we have
\begin{align}
 F'_1(\bubar(\bx)\cdot \be/\sqrt{s}) - F'_1(\bvbar(\bx)\cdot \be/\sqrt{s}) &= \int_0^1 \dparder{F'_1(\br(l,\bx) \cdot \be/\sqrt{s})}{l} dl \\
 &= \int_0^1 \dparder{F'_1(\br \cdot \be/\sqrt{s})}{\br} \vert_{\br= \br(l,\bx)} \cdot \dparder{\br(l,\bx)}{l} dl. \label{eq:estimate F prime bx}
\end{align}

Similarly, we have.
\begin{align}
 F'_1(\bubar(\by)\cdot \be/\sqrt{s}) - F'_1(\bvbar(\by)\cdot \be/\sqrt{s}) &= \int_0^1 \dparder{F'_1(\br \cdot \be/\sqrt{s})}{\br} \vert_{\br= \br(l,\by)} \cdot \dparder{\br(l,\by)}{l} dl. \label{eq:estimate F prime by}
\end{align}
Note that 
\begin{align}\label{eq:F prime wrt br}
\dparder{F'_1(\br \cdot \be/\sqrt{s})}{\br} \vert_{\br= \br(l,\by)} &= F''_1(\br(l,\bx) \cdot \be/\sqrt{s}) \dfrac{\be}{\sqrt{s}}.
\end{align}

Combining \autoref{eq:estimate F prime bx}, \autoref{eq:estimate F prime by}, and \autoref{eq:F prime wrt br}, gives
\begin{align*}
H &= \dfrac{1}{\abs{\bx - \by}^\gamma} \vline \int_0^1 \left(\tilde{\omega}(\bx,\bxi) F''_1(\br(l,\bx) \cdot \be/\sqrt{s}) (\bubar(\bx) - \bvbar(\bx)) \right. \notag \\
& \qquad \qquad \left. - \tilde{\omega}(\by,\bxi)F''_1(\br(l,\by) \cdot \be/\sqrt{s}) (\bubar(\by) - \bvbar(\by)) \right) \cdot \dfrac{\be}{\sqrt{s}} dl \vline \notag \\
&\leq \dfrac{1}{\abs{\bx - \by}^\gamma} \dfrac{1}{\sqrt{s}} | \int_0^1 | \tilde{\omega}(\bx,\bxi)F''_1(\br(l,\bx) \cdot \be/\sqrt{s}) (\bubar(\bx) - \bvbar(\bx))  \notag \\
& \qquad \qquad  - \tilde{\omega}(\by,\bxi)F''_1(\br(l,\by) \cdot \be/\sqrt{s}) (\bubar(\by) - \bvbar(\by)) | dl |. \notag \\
\end{align*}
Adding and subtracting $\tilde{\omega}(\bx,\bxi)F''_1(\br(l,\bx)\cdot\be/\sqrt{s})(\bubar(\by) - \bvbar(\by))$, and noting $0\leq \tilde{\omega}(\bx,\bxi)\leq 1$ gives
\begin{align*}
H &\leq \dfrac{1}{\abs{\bx - \by}^\gamma} \dfrac{1}{\sqrt{s}} | \int_0^1 |F''_1(\br(l,\bx)\cdot\be/\sqrt{s})| \abs{\bubar(\bx) - \bvbar(\bx) - \bubar(\by) + \bvbar(\by)} dl| \notag \\
&\quad + \dfrac{1}{\abs{\bx - \by}^\gamma} \dfrac{1}{\sqrt{s}}  \int_0^1 |(\tilde{\omega}(\bx,\bxi)F''_1(\br(l,\bx)\cdot\be/\sqrt{s}) - \tilde{\omega}(\by,\bxi)F''_1(\br(l,\by)\cdot\be/\sqrt{s}))| \notag\\
&\quad \quad \quad \quad \quad \quad \quad \quad\times \abs{\bubar(\by) - \bvbar(\by)} dl . \notag \\
&=: H_1 + H_2.
\end{align*}

Estimating $H_1$ first. Note that $|F''_1(r)|\leq C_2$. Since $\bu,\bv \in \Cholderz{\perdrd}$, it is easily seen that
\begin{align*}
\dfrac{| \bubar(\bx) - \bvbar(\bx) - \bubar(\by) + \bvbar(\by) |}{\abs{\bx - \by}^\gamma} &\leq 2 \Choldernorm{\bu - \bv}{\perdrd}.
\end{align*}
Therefore, we have
\begin{align}
H_1 &\leq \dfrac{2 C_2}{\sqrt{s}} \Choldernorm{\bu - \bv}{\perdrd}. \label{eq:estimate H one}
\end{align}

We now estimate $H_2$. We add and subtract $ \tilde{\omega}(\bx,\bxi)F''_1(\br(l,\by)\cdot\be/\sqrt{s}))$ in $H_2$ to get
\begin{align*}
H_2 \leq H_3+H_4, 
\end{align*}
where
\begin{align*}
H_3 = & \dfrac{1}{\abs{\bx - \by}^\gamma} \dfrac{1}{\sqrt{s}}  \int_0^1 |(F''_1(\br(l,\bx)\cdot\be/\sqrt{s}) - F''_1(\br(l,\by)\cdot\be/\sqrt{s}))| \abs{\bubar(\by) - \bvbar(\by)} dl,\notag
\end{align*}
and
\begin{align*}
H_4 &= \dfrac{1}{\abs{\bx - \by}^\gamma} \dfrac{1}{\sqrt{s}}  \int_0^1 |(\tilde{\omega}(\bx,\bxi) - \tilde{\omega}(\by,\bxi)|F''_1(\br(l,\by)\cdot\be/\sqrt{s}))| \abs{\bubar(\by) - \bvbar(\by)} dl . \notag
\end{align*}
Now we estimate $H_3$. Since $|F'''_1(r)| \leq C_3$, see \autoref{C3}, we have
\begin{align}
&\dfrac{1}{\abs{\bx - \by}^\gamma}|F''_1(\br(l,\bx)\cdot\be/\sqrt{s}) - F''_1(\br(l,\by)\cdot\be/\sqrt{s})| \notag \\
&\leq \dfrac{1}{\abs{\bx - \by}^\gamma} \sup_r \abs{F'''(r)} \dfrac{\abs{\br(l,\bx)\cdot \be - \br(l,\by) \cdot \be} }{\sqrt{s}} \notag \\
&\leq \dfrac{C_3}{\sqrt{s}} \dfrac{\abs{\br(l,\bx) - \br(l,\by)}}{\abs{\bx - \by}^\gamma}  \notag \\
&= \dfrac{C_3}{\sqrt{s}} \left( \dfrac{\abs{1-l} \abs{\bvbar(\bx) - \bvbar(\by)}}{\abs{\bx - \by}^\gamma} + \dfrac{\abs{l} \abs{\bubar(\bx) - \bubar(\by)}}{\abs{\bx - \by}^\gamma} \right) \notag \\
&\leq \dfrac{C_3}{\sqrt{s}} \left( \dfrac{\abs{\bvbar(\bx) - \bvbar(\by)}}{\abs{\bx - \by}^\gamma} + \dfrac{ \abs{\bubar(\bx) - \bubar(\by)}}{\abs{\bx - \by}^\gamma} \right). \label{eq:estimate H two part one}
\end{align}
Where we have used the fact that $\abs{1-l} \leq 1, \abs{l} \leq 1$, as $l\in [0,1]$. Also, note that
\begin{align*}
\dfrac{\abs{\bubar(\bx) - \bubar(\by)}}{\abs{\bx - \by}^{\gamma}} &\leq 2 \Choldernorm{\bu}{\perdrd} \\
\dfrac{\abs{\bvbar(\bx) - \bvbar(\by)}}{\abs{\bx - \by}^{\gamma}} &\leq 2 \Choldernorm{\bv}{\perdrd} \\
\abs{\bubar(\by) - \bvbar(\by)} &\leq s^\gamma \Choldernorm{\bu - \bv}{\perdrd}.
\end{align*}
We combine above estimates with \autoref{eq:estimate H two part one}, to get
\begin{align}
H_3&\leq \dfrac{1}{\sqrt{s}} \dfrac{C_3}{\sqrt{s}} \left( \Choldernorm{\bu}{\perdrd} + \Choldernorm{\bv}{\perdrd} \right) s^\gamma \Choldernorm{\bu - \bv}{\perdrd} \notag \\
&= \dfrac{C_3}{s^{1-\gamma}} \left( \Choldernorm{\bu}{\perdrd} + \Choldernorm{\bv}{\perdrd} \right) \Choldernorm{\bu - \bv}{\perdrd}.\label{eq:estimate H two}
\end{align}

Next we estimate $H_4$. Here we add and subtract $\omega(\by)\omega(\bx+\epsilon\bxi)$ to get
\begin{align*}
H_4 &= \dfrac{1}{\abs{\bx - \by}^\gamma} \dfrac{1}{\sqrt{s}}  \int_0^1 |(\omega(\bx,\bx+\epsilon\bxi)(\omega(\bx)-\omega(\by)) +\omega(\by)(\omega(\bx+\epsilon\bxi)-\omega(\by+\epsilon\bxi))\notag\\
&\quad \quad \quad \quad \quad \quad \quad \quad\times|F''_1(\br(l,\by)\cdot\be/\sqrt{s}))| \abs{\bubar(\by) - \bvbar(\by)} dl . \notag
\end{align*}
Recalling that $\omega$ belongs to $\Cholderz{\perdrd}$ and in view of the previous estimates a straight forward calculation gives
\begin{align}
H_4\leq \dfrac{4C_2}{s^{1/2-\gamma}}\Choldernorm{\omega}{\perdrd}\Choldernorm{\bu - \bv}{\perdrd}.\label{eq:estimate H three} 
\end{align}

Combining \autoref{eq:estimate H one}, \autoref{eq:estimate H two}, and \autoref{eq:estimate H three}  gives
\begin{align*}
H &\leq \left( \dfrac{2C_2}{\sqrt{s}} + \dfrac{4C_2}{s^{1/2-\gamma}}\Choldernorm{\omega}{\perdrd}+  \right.\notag\\
&\left.+\dfrac{C_3}{s^{1-\gamma}} \left( \Choldernorm{\bu}{\perdrd} + \Choldernorm{\bv}{\perdrd} \right) \right) \Choldernorm{\bu - \bv}{\perdrd}.
\end{align*}

Substituting $H$ in \autoref{eq:H def}, gives
\begin{align}
&\dfrac{1}{\abs{\bx - \by}^{\gamma}} \abs{(-\del PD^{\epsilon}(\bu) + \del PD^{\epsilon}(\bv))(\bx) - (-\del PD^{\epsilon}(\bu) + \del PD^{\epsilon}(\bv))(\by)} \notag \\
&\leq | \dfrac{2}{\epsilon\omega_d} \int_{H_1(\bzero)} J(\abs{\bxi}) \dfrac{1}{\sqrt{s}} H  d\bxi | \notag \\
&\leq \left( \dfrac{4C_2 \bar{J}_1}{\epsilon^2} +\dfrac{4C_2\bar{J}_{1-\gamma}}{\epsilon^{2-\gamma}} \Choldernorm{\omega}{\perdrd}\right.\notag\\
&\left.+ \dfrac{2C_3 \bar{J}_{3/2 - \gamma}}{\epsilon^{2+1/2 - \gamma}} \left( \Choldernorm{\bu}{\perdrd} + \Choldernorm{\bv}{\perdrd} \right)  \right) \Choldernorm{\bu - \bv}{\perdrd}.
\label{eq:estimate lipschitz part 2}
\end{align}
We combine \autoref{eq:lipschitz norm of per force}, \autoref{eq:estimate lipschitz part 1}, and \autoref{eq:estimate lipschitz part 2}, and get
\begin{align}
&\Choldernorm{-\del PD^{\epsilon}(\bu) - (-\del PD^{\epsilon}(\bv))}{} \notag \\
&\leq \left( \dfrac{4C_2 \bar{J}_1}{\epsilon^2} +\dfrac{2C_2 \bar{J}_{1-\gamma}}{\epsilon^{2-\gamma}} (1+\Choldernorm{\omega}{}) + \dfrac{2C_3 \bar{J}_{3/2 - \gamma}}{\epsilon^{2+1/2 - \gamma}} \left( \Choldernorm{\bu}{} + \Choldernorm{\bv}{} \right)   \right) \Choldernorm{\bu - \bv}{} \notag \\
&\leq \dfrac{\bar{C}_1 + \bar{C}_2\Choldernorm{\omega}{} +\bar{C}_3(\Choldernorm{\bu}{} + \Choldernorm{\bv}{})}{\epsilon^{2+\alpha(\gamma)}} \Choldernorm{\bu - \bv}{} \label{eq:del pd lipschitz}
\end{align}
where we introduce new constants $\bar{C}_1, \bar{C}_2, \bar{C}_3$. We let $\alpha(\gamma) = 0$, if $\gamma \geq 1/2$, and $\alpha(\gamma) = 1/2 - \gamma$, if $\gamma \leq 1/2$. One can easily verify that, for all $\gamma \in (0,1]$ and $0< \epsilon \leq 1$,
\begin{align*}
\max \leftcr\dfrac{1}{\epsilon^2}, \dfrac{1}{\epsilon^{2+1/2 - \gamma}}, \dfrac{1}{\epsilon^{2-\gamma}} \rightcr \leq \dfrac{1}{\epsilon^{2+ \alpha(\gamma)}}
\end{align*}

To complete the proof, we combine \autoref{eq:del pd lipschitz} and \autoref{eq:norm of F in X}, and get 
\begin{align*}
\normX{F^\epsilon(y,t) - F^\epsilon(z,t)}{X} &\leq \dfrac{L_1 + L_2 (\Choldernorm{\omega}{}+\normX{y}{X} + \normX{z}{X})}{\epsilon^{2+\alpha(\gamma)}} \normX{y-z}{X}.
\end{align*}
This proves the Lipschitz continuity of $F^\epsilon(y,t)$ on any bounded subset of $X$. 
The bound on $F^\epsilon(y,t)$, see \autoref{eq:bound on F}, follows easily from \autoref{eq:del pd simplified}. This completes the proof of \autoref{prop:lipschitz}. 

\subsection{Existence of solution in H\"{o}lder space}
\label{section: existence}
In this section, we prove \autoref{thm:existence over finite time domain}. We begin by proving a local existence theorem. We then show that the local solution can be continued uniquely in time to recover \autoref{thm:existence over finite time domain}.

The existence and uniqueness of local solutions is stated in the following theorem.

{\vskip 2mm}
\begin{theorem}\label{thm:local existence}
\textbf{Local existence and uniqueness} \\
Given $X= \Cholderz{\perdrd}\times\Cholderz{\perdrd}$, $\bb(t)\in C^{0,\gamma}_0(D;\mathbb{R}^d)$, and  initial data $x_0=(\bu_0,\bv_0)\in X$. We suppose that $\bb(t)$ is continuous in time over some time interval $I_0=(-T,T)$ and satisfies $\sup_{t\in I_0} \Choldernorm{\bb(t)}{\perdrd} < \infty$.  Then, there exists a time interval $I'=(-T',T')\subset I_0$ and unique solution $y =(y^1,y^2) $ such that $y\in C^1(I';X)$ and

\begin{equation}
y(t)=x_0+\int_0^tF^\epsilon(y(\tau),\tau)\,d\tau,\hbox{  for $t\in I'$}
\label{8loc}
\end{equation}

or equivalently

\begin{equation*}
y'(t)=F^\epsilon(y(t),t),\hbox{with    $y(0)=x_0$},\hbox{  for $t\in I'$}
\label{11loc}
\end{equation*}
where $y(t)$ and $y'(t)$ are Lipschitz continuous in time for $t\in I'\subset I_0$.

\end{theorem}
{\vskip 2mm}

\begin{figure} 
\centering
\begin{tikzpicture}[xscale=1.0,yscale=1.0]
\draw [] (0.0,0.0) circle [radius=2.9];
\draw [] (1.0,2.0) circle [radius=0.6];
\draw [->,thick] (1.0,2.0) -- (1.38729,2.38729);
\node [right] at (0.8,1.8) {$x_0$};
\node [below] at (-0.2,1.82) {$B(x_0,r)$};
\node [below] at (0.0,0.0) {$0$};
\draw [->,thick] (0,0) -- (2.9,0.0);
\node [below] at (1.45,0.0) {$R$};
\node [below] at (-2.5,-2.3) {$B(0,R)$};
\end{tikzpicture} 
\caption{Geometry.}
 \label{figurenested}
\end{figure}
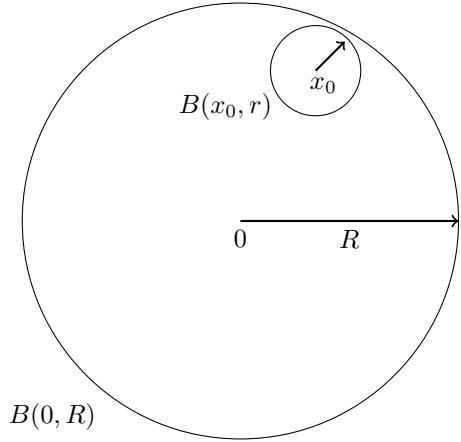

To prove \autoref{thm:local existence}, we proceed as follows. 
We write $y(t)=(y^1(t),y^2(t))$ and $||y||_X=||y^1(t)||_{\Cholder{}}+||y^2(t)||_{\Cholder{}}$. 
Define the ball $B(0,R)= \{y\in X:\, ||y||_X<R\}$ and choose $R>||x_0||_X$. Let $r = R - \normX{x_0}{X}$  and we consider the ball $B(x_0,r)$ defined by
\begin{equation}
B(x_0,r)=\{y\in X:\,||y-x_0||_X<r\}\subset B(0,R),
\label{balls}
\end{equation}
see figure \autoref{figurenested}. 

To recover the existence and uniqueness we introduce the transformation
\begin{equation*}
S_{x_0}(y)(t)=x_0+\int_0^tF^\epsilon(y(\tau),\tau)\,d\tau.
\label{0}
\end{equation*}
Introduce $0<T'<T$ and the associated set $Y(T')$ of H\"older continuous functions taking values in $B(x_0,r)$ for  $I'=(-T',T')\subset I_0=(-T,T)$.
The goal is to find appropriate interval $I'=(-T',T')$ for which $S_{x_0}$ maps into the corresponding set $Y(T')$. Writing out the transformation with $y(t)\in Y(T')$ gives
\begin{eqnarray}
&&S_{x_0}^1(y)(t)=x_0^1+\int_0^t y^2(\tau)\,d\tau\label{1}\\
&&S_{x_0}^2(y)(t)=x_0^2+\int_0^t(-\nabla PD^\epsilon(y^1(\tau))+\bb(\tau))\,d\tau,\label{2}
\end{eqnarray}
and there is a positive constant $K=C/\epsilon^{2+\alpha(\gamma)}$, see 
 \autoref{eq:bound on F}, independent of $y^1(t)$, for $-T'<t<T'$, such that estimation in \autoref{2} gives
\begin{eqnarray}
||S_{x_0}^2(y)(t)-x_0^2||_{\Cholder{}}\leq (K(1+\frac{1}{\epsilon^\gamma}+\sup_{t\in(-T',T')}||y^1(t)||_{\Cholder{}})+\sup_{t\in(-T,T)}||\bb(t)||_{\Cholder{}})T'\nonumber\\
\label{3}
\end{eqnarray}
and from \autoref{1}
\begin{eqnarray}
||S_{x_0}^1(y)(t)-x_0^1||_{\Cholder{}}\leq\sup_{t\in(-T',T')}||y^2(t)||_{\Cholder{}}T'\label{4}.
\end{eqnarray}

We write $b=\sup_{t\in I_0}||\bb(t)||_{\Cholder{}}$ and adding \autoref{3} and \autoref{4} gives the upper bound
\begin{equation}
||S_{x_0}(y)(t)-x_0||_X\leq (K(1+\frac{1}{\epsilon^\gamma}+\sup_{t\in(-T',T')}||y(t)||_X)+b)T'.
\label{5}
\end{equation}Since $B(x_0,r)\subset B(0,R)$, see \autoref{balls}, we make the choice $T'$ so that
\begin{equation}
||S_{x_0}(y)(t)-x_0||_X\leq ((K(1+\frac{1}{\epsilon^\gamma}+R)+b)T'<r=R-||x_0||_X.
\label{5point5}
\end{equation}
For this choice we see that
\begin{equation}
T'<\theta(R)=\frac{R-||x_0||_X}{K(R+1+\frac{1}{\epsilon^\gamma})+b}.
\label{6}
\end{equation}
Now it is easily seen that $\theta(R)$ is increasing with $R>0$ and 
\begin{equation}
\lim_{R\rightarrow\infty}\theta(R)=\frac{1}{K}.
\label{7}
\end{equation}
So given $R$ and $||x_0||_X$ we choose $T'$ according to
\begin{equation}
\frac{\theta(R)}{2}<T'< \theta(R),
\label{localchoiceofT}
\end{equation}
and set $I'=(-T',T')$. We have found the appropriate time domain $I'$ such that the transformation $S_{x_0}(y)(t)$ as defined in \autoref{0} maps $Y(T')$ into itself. We now proceed using standard arguments, see e.g. [\cite{MA-Driver}, Theorem 6.10], to complete the proof of existence and uniqueness of solution for given initial data $x_0$ over the interval $I'= (-T', T')$.

We now prove \autoref{thm:existence over finite time domain}. From the proof of \autoref{thm:local existence}  above, we see that a unique local solution exists over a time domain $(-T',T')$ with $\frac{\theta(R)}{2}<T'$. Since $\theta(R)\nearrow1/K$ as $R\nearrow \infty$ we can fix a tolerance $\eta>0$ so that $[(1/2K)-\eta]>0$. Then given any initial condition with bounded H\"older norm and $b=\sup_{t\in[-T,T)}||\bb(t)||_{\Cholder{}}$ we can choose $R$ sufficiently large so that $||x_0||_X<R$ and   $0<(1/2K))-\eta<T'$. Thus we can always find local solutions for time intervals $(-T',T')$ for $T'$ larger than $[(1/2K)-\eta]>0$. Therefore we apply the local existence and uniqueness result to uniquely continue local solutions up to an arbitrary time interval $(-T,T)$.

\section{Limit behavior of H\"older solutions in the limit of vanishing nonlocality}
\label{s:discussion}
In this section, we consider the behavior of bounded H\"older continuous solutions as the peridynamic horizon tends to zero. We find that the solutions converge to a limiting sharp fracture evolution with bounded Griffiths fracture energy and satisfy the linear elastic wave equation away from the fracture set. We look at a subset of H\"older solutions that are differentiable in the spatial variables to show that sharp fracture evolutions can be approached by spatially smooth evolutions in the limit of vanishing non locality. As $\epsilon$ approaches zero derivatives can become large but must localize to surfaces across which the limiting evolution jumps.

We consider a sequence of peridynamic horizons $\epsilon_k=1/k$, $k=1,\ldots$ and the associated H\"older continuous solutions $\bu^{\epsilon_k}(t,\bx)$ of the peridynamic initial value problem \autoref{eq:per equation},  \autoref{eq:per bc}, and \autoref{eq:per initialvalues}. We assume that the initial conditions $\bu_0^{\epsilon_k},\bv_0^{\epsilon_k}$ have uniformly bounded peridynamic energy and mean square initial velocity given by
\begin{equation*}
\sup_{\epsilon_k}PD^{\epsilon_k}(\bu^{\epsilon_k}_0)<\infty\hbox{  and  }\sup_{\epsilon_k}||\bv^{\epsilon_k}_0||_{L^2(D;\mathbb{R}^d)}<\infty.
\end{equation*}
Moreover we suppose that $\bu_0^{\epsilon_k},\bv_0^{\epsilon_k}$  are differentiable on $D$ and that they converge in $L^2(D;\mathbb{R})$ to $\bu_0^{0},\bv_0^{0}$ with bounded Griffith free energy given by
\begin{eqnarray*}
\int_{D}\,2\mu |\mathcal{E} \bu_0^0|^2+\lambda |{\rm div}\,\bu_0^0|^2\,dx+\mathcal{G}_c\mathcal{H}^{d-1}(J_{\bu_0^0})\leq C < \infty, 
\end{eqnarray*}
where $J_{\bu_0^0}$ denotes an initial fracture surface given by the jumps in the initial deformation $\bu_0^0$ and $\mathcal{H}^{2}(J_{\bu^0(t)})$ is its $2$ dimensional Hausdorff measure of the jump set. Here $\mathcal{E} \bu^0_0$ is the elastic strain and ${\rm div}\,\bu^0_0=Tr(\mathcal{E} \bu^0_0)$. The constants $\mu$, $\lambda$ are given by the explicit formulas
\begin{eqnarray*}
\hbox{ and }& \mu=\lambda=\frac{1}{5} f'(0)\int_{0}^1r^dJ(r)dr, \hbox{ $d=2,3$}
\end{eqnarray*}
and
\begin{eqnarray*}
\mathcal{G}_c=\frac{3}{2}\, f_\infty \int_{0}^1r^dJ(r)dr, \hbox{ $d=2,3$},
\end{eqnarray*}
where $f'(0)$ and $f_\infty$ are defined by \autoref{eq:per asymptote}. Here $\mu=\lambda$ and is a consequence of the central force model used in cohesive dynamics. 
Last we suppose as in \cite{CMPer-Lipton} that the solutions are uniformly bounded, i.e.,
\begin{equation*}
\sup_{\epsilon_k}\sup_{[0,T]}||\bu^{\epsilon_k}(t)||_{L^\infty(D;\mathbb{R}^d)}<\infty,
\end{equation*}

The H\"older solutions $\bu^{\epsilon_k}(t,\bx)$ naturally belong to $L^2(D;\mathbb{R}^d)$ for all $t\in[0,T]$ and we can directly apply the Gr\"onwall inequality (equation (6.9) of \cite{CMPer-Lipton}) together with Theorems 6.2  and 6.4 of \cite{CMPer-Lipton} to conclude similar to Theorems 5.1 and 5.2  of \cite{CMPer-Lipton} that there is at least one ``cluster point'' $\bu^{0}(t,\bx)$ belonging to $C([0,T];L^2(D;\mathbb{R}^d))$ and subsequence, also denoted by $\bu^{\epsilon_k}(t,\bx)$ for which 
\begin{eqnarray*}
\lim_{\epsilon_k\rightarrow 0}\max_{0\leq t\leq T}\left\{\Vert \bu^{\epsilon_k}(t)-\bu^0(t)\Vert_{L^2(D;\mathbb{R}^d)}\right\}=0.
\label{eq:per unifconvg}
\end{eqnarray*}
Moreover it follows from \cite{CMPer-Lipton} that the limit evolution $\bu^0(t,\bx)$ has a weak derivative $\bu_t^0(t,\bx)$ belonging to $L^2([0,T]\times D;\mathbb{R}^{d})$. For each time $t\in[0,T]$ we can apply methods outlined in \cite{CMPer-Lipton} to find that the cluster point $\bu^0(t,\bx)$  is a special function of bounded deformation (see, \cite{Ambrosio}, \cite{Bellettini}) and has bounded linear elastic fracture energy given by
\begin{eqnarray*}
\int_{D}\,2\mu |\mathcal{E} \bu^0(t)|^2+\lambda |{\rm div}\,\bu^0(t)|^2\,dx+\mathcal{G}_c\mathcal{H}^{2}(J_{\bu^0(t)})\leq C,
\label{LEFMbound}
\end{eqnarray*}
for $0\leq t\leq T$ where $J_{\bu^0(t)}$ denotes the evolving fracture surface 
The deformation - crack set pair $(\bu^0(t),J_{\bu^0(t)})$ records the brittle fracture evolution of the limit dynamics.

Arguments identical to \cite{CMPer-Lipton} show that away from sets where $|S(\by,\bx;\bu^{\epsilon_k})|>S_c$ the limit $\bu^0$ satisfies the linear elastic wave equation. This is stated as follows:
Fix $\delta>0$ and for   $\epsilon_k<\delta$ and $0\leq t\leq T$ consider the open set $D'\subset D$  for which points $\bx$ in $D'$  and $\by$ for which $|\by-\bx|<\epsilon_k$ satisfy,
\begin{eqnarray*}
|S(\by,\bx;\bu^{\epsilon_k}(t))|<{S}_c(\by,\bx).
\label{eq: per quiecent}
\end{eqnarray*}
Then the limit evolution $\bu^0(t,\bx)$  evolves elastodynamically  on $D'$ and is governed by the balance of linear momentum expressed by the Navier Lam\'e equations on the domain $[0,T]\times D'$ given by
\begin{eqnarray*}
 \bu^0_{tt}(t)= {\rm div}\bolds{\sigma}(t)+\bb(t), \hbox{on $[0,T]\times D'$},
\label{waveequationn}
\end{eqnarray*}
where the stress tensor $\bolds{\sigma}$ is given by,
\begin{eqnarray*}
\bolds{\sigma} =\lambda I_d Tr(\mathcal{E}\,\bu^0)+2\mu \mathcal{E}\bu^0,
\label{stress}
\end{eqnarray*}
where $I_d$ is the identity on $\mathbb{R}^d$ and $Tr(\mathcal{E}\,\bu^0)$ is the trace of the strain. 
Here the second derivative $\bu_{tt}^0$ is the time derivative in the sense of distributions of $\bu^0_t$ and ${\rm div}\bolds{\sigma}$ is the divergence of the stress tensor $\bolds{\sigma}$ in the distributional sense. This shows that sharp fracture evolutions can be approached by spatially smooth evolutions in the limit of vanishing non locality.

\section{Conclusions}
\label{s:conclusions}
In this article, we have presented a numerical analysis for class of nonlinear nonlocal  peridynamic models. We have shown that the convergence rate applies, even when the fields do not have well-defined spatial derivatives. We treat both the forward Euler scheme as well as the general implicit single step method. 
The convergence rate is found to be the same for both schemes and is given by $C(\Delta t+h^\gamma/\epsilon^2)$. Here the constant $C$ depends on $\epsilon$ and H\"older and $L^2$ norm of the solution and its time derivatives.
The Lipschitz property of the nonlocal, nonlinear force together with boundedness of the nonlocal kernel plays an important role. It ensures that the error in the nonlocal force remains bounded when replacing the exact solution with its approximation. This, in turn, implies that even in the presence of mechanical instabilities the global approximation error remains controlled by the local truncation error in space and time. 

Taking $\gamma=1$, a straight forward estimate of \autoref{eq:const Cs} using \autoref{eq:bound on F} gives to leading order
\begin{align}\label{eq: final est initial}
\sup_{0\leq k \leq T/\Delta t} E^k\leq \left [C_1 \Delta t C_t + C_2 h\sup_{0<t<T}\Vert u\Vert_{C^{0,1}(D;\mathbb{R}^3)}\right ],
\end{align}
where $C_t$ is independent of $\epsilon$ and depends explicitly on the $L^2$ norms of time derivatives of the solution see \autoref{eq:const Ct} and 
\begin{align*}
C_1=\exp \left[T (1 + 6\bar{C}/\epsilon^2) \right]T,
\end{align*}
\begin{align*}
C_2=\exp \left[T (1 + 6\bar{C}/\epsilon^2) \right]T\left(1+\sqrt{3}\bar{C}(1+\frac{1}{\epsilon})+\frac{4\sqrt{3}\bar{C}}{\epsilon^2}\right).
\end{align*}
It is evident that the exponential factor could be large. However we can choose times $T$  for which the effects of the exponential factor can be diminished and $C_1$ and $C_2$ are not too large.  To fix ideas consider a $1$ cubic meter sample and a corresponding 1400 meter per second shear wave speed. This wave speed is characteristic of plexiglass.
Then the time for a shear wave to traverse the sample is $718$ $\mu$-seconds. This is the characteristic time $T^\ast$ and a fracture experiment can last a few hundred $\mu$-seconds.  The actual time in $\mu$-seconds of  a fracture simulation is given by $TT^*$ where $T$ is the non-dimensional simulation time. The dimensionless constant $\bar{C}$ is $1.19$ and we take $\epsilon=1/10$ and dimensionless body force unity. For a simulation cycle of length $TT^\ast=1.5\mu$-seconds the constants $C_1$ and  $C_2$ in \autoref{eq: final est initial} are $0.0193$ and $7.976$ respectively. The solution after cycle time $T$ can be used as initial conditions for a subsequent run and the process can be iterated. Unfortunately these estimates predict a total simulation time of $15\mu$-second before the relative error becomes greater than one even for a large number of spatial degrees of freedom. We point out that because the constants in the a-priori bound are necessarily pessimistic the predicted simulation time is an order of magnitude below what is seen in experiment. Future work will focus on a-posteriori estimates for simulations and adaptive implementations of the finite difference scheme.

In conclusion the analysis shows that the method is stable and one can control the error by choosing the time step and spatial discretization sufficiently small.  However errors do accumulate with time steps and this limits the time interval for simulation. We have identified local perturbations for which the error accumulates with time step for the implicit Euler method. These unstable local perturbations  correspond to regions for which a preponderance of bonds are in the softening regime.

\section*{Acknowledgements}
RL would like to acknowledge the support and kind hospitality of the Hausdorff Institute for Mathematics in Bonn during the trimester program on multiscale problems.


\newcommand{\noopsort}[1]{}

\end{document}